\documentstyle{article}

\def\kon{\wedge}
\def\str{\rightarrow}
\def\rts{\leftarrow}
\def\mj{\mbox{\bf 1}}
\def\cm{\mbox{$SMI(\cal P)$}}
\def\cmst{$\cm^{st}$}
\def\ca{\mbox{$\cal A$}}
\def\cirk{\,{\raisebox{.3ex}{\tiny $\circ$}}\,}
\def\mN{\mbox{$\mathbf{N}$}}
\def\prop#1#2{\vspace{2ex} \noindent{\sc #1.} {\it #2} \par \vspace{2ex}}
\def\dkz{\noindent{\sc Proof. }}
\def\qed{\hfill $\dashv$}
\def\lz{|\![}
\def\dz{]\!|}
\def\pl{\!+\!}
\def\mn{\!-\!}
\def\strt{\stackrel{\textbf{.}\,}{\rightarrow}}

\begin{document}

\title{Symmetric bimonoidal intermuting categories and $\omega\times\omega$ reduced bar constructions}
\author{\small {\sc Zoran Petri\' c}\footnote{Mathematical Institute, SANU, Knez Mihailova 36, p.f.\ 367, 11001 Belgrade,
Serbia, email: zpetric@mi.sanu.ac.rs}\hspace{2em}{\sc Todd
Trimble}\footnote{8 Iris Lane Redding, CT 06896 USA, email:
topological.musings@gmail.com}}
\date{}
\maketitle

\vspace{-3ex}

\begin{abstract}
\noindent A new, self-contained, proof of a coherence result for
categories equipped with two symmetric monoidal structures bridged
by a natural transformation is given. It is shown that this
coherence result is sufficient for $\omega\times\omega$-indexed
family of iterated reduced bar constructions based on such a
category.

\end{abstract}

\vspace{.3cm}


\vspace{.5ex}

\noindent {\small {\it Keywords$\,$}: coherence, bar construction,
symmetric monoidal categories, infinite loop spaces}

\vspace{.5ex}

\noindent {\small {\it Acknowledgements$\,$}: This work was
supported by a project of the Ministry of Science of Serbia
(ON174026). }

\section{\large\bf Introduction}
This paper presents a reduced bar construction which is usually
the initial part of the results connecting various monoidal
categories with 1-fold, 2-fold, $n$-fold and infinite loop spaces
(see \cite{S74}, \cite{T79}, \cite{BFSV} and references therein).
By the reduced bar construction we mean a construction of a
simplicial object based on a monoid in a category whose monoidal
structure is given by finite products (exactly the same as the
notion used in \cite{T79}), which in particular, for a special
monoid in the category \emph{Cat}, may be iterated in order to
obtain a lax functor from an arbitrary power of the opposite of
topologist's simplicial category to \emph{Cat}. This construction
is based here on a category equipped with two symmetric monoidal
structures, given by the tensors $\vee$ and $\kon$, and the units
$\bot$ and $\top$. These two structures are bridged by a natural
transformation, called \emph{intermutation} in \cite{DPInt}, given
by the family of arrows
\[
(A\kon B)\vee(C\kon D)\str(A\vee C)\kon(B\vee D).
\]
Such categories appeared under the name \emph{symmetric bimonoidal
intermuting categories} in \cite{DPInt}. As a concrete example of
a symmetric bimonoidal intermuting category one can take any
category with all finite coproducts and all finite products in
which product of initial objects is initial and coproduct of
terminal objects is terminal (see \cite{DPInt}, Sections 13 and
15).

We will not go further in the procedure of delooping, which is
very well traced by the work of Thomason in \cite{T79}. This
procedure leads to an $\omega\times\omega$-indexed family of
deloopings of the classifying space of a symmetric bimonoidal
intermuting category. (According to this, one can make the
following hierarchy of infinite loop spaces; simply symmetric
monoidal structure corresponds to an infinite loop space with
$\omega$-indexed family of deloopings, double symmetric monoidal
structure without intermutation corresponds to an infinite loop
space with $\omega+\omega$-indexed family of deloopings, double
symmetric monoidal structure with intermutation corresponds to an
infinite loop space with $\omega\times\omega$-indexed family of
deloopings, etc.)

A definition of $n$-fold monoidal category is usually inductive
and it starts with pseudomonoids (or monoids) in the 2-category
\emph{Cat} whose monoidal structure is given by 2-products. Then
one makes a choice what to consider to be the morphisms between
monoidal (or strict monoidal) categories, i.e.\ how strict they
should preserve the monoidal structure. This leads to a 2-category
$Mon(Cat)$, again with 2-products. A pseudomonoid (or a monoid) in
such a category is a (strict) 2-fold monoidal category and if we
iterate the above with the same degree of strictness, we obtain
one possible notion of $n$-fold monoidal category.

In \cite{JS93}, Joyal and Street deal with such a concept having
in its basis the 2-category $Mon_{ps}(Cat)$, which is the
2-monoidal category of monoidal categories, ``pseudo'' or strong
monoidal functors, and monoidal transformations. They have shown
that such a degree of strictness leads to a sequence of categorial
structures starting with monoidal categories, then we have braided
monoidal categories as 2-fold monoidal categories and symmetric
monoidal categories as n-fold monoidal categories for $n\geq 3$.
In terms of loop spaces, these categorial structures model (up to
group completion) 1-fold loop spaces, 2-fold loop spaces, and
infinite loop spaces.

In \cite{BFSV}, Balteanu, Fiedorowicz, Schw\" anzl and Vogt
considered a variant of $Mon(Cat)$ in which the interchange
between multiplicative structures need not be invertible. This was
an important advance leading to a definition of $n$-fold monoidal
categories without stabilization at $n=3$. However, they did not
go far enough by similarly laxifying the appropriate interchanges
for units, which were treated in their work as strict as possible.

Let $Mon_{lax}(Cat)$ be the variant of $Mon(Cat)$ in which the
interchange between multiplicative structures and interchange
between units need not be invertible, i.e.\ a 2-monoidal category
of monoidal categories, lax monoidal functors, and monoidal
transformations. The possibility of defining $n$-fold monoidal
structures with respect to such a basis is much less explored
perhaps because of difficulties in proving corresponding coherence
results.

Here we deal with a categorial structure obtained by an analogous
iteration. We start with $SyMon_{lax}(Cat)$, a 2-monoidal category
of symmetric monoidal categories, lax symmetric monoidal functors
and monoidal transformations. At the next level we have a
2-monoidal category $SyMon^2_{lax}(Cat)$ of pseudocommutative
pseudomonoids in $SyMon_{lax}(Cat)$, i.e.\ 2-fold symmetric
monoidal categories. (By iterating this, one could define a notion
of $n$-fold symmetric monoidal category.) Our symmetric bimonoidal
intermuting categories are objects of $SyMon_{lax}^{2}(Cat)$ for
which we require some structural constraints to be invertible.
However, we can't find many convincing examples of symmetric
bimonoidal intermuting categories. (Here we deal with a
syntactically built one and its variants.) We still have no
general coherence result that provides the desired reduced bar
construction based on an arbitrary object of
$SyMon_{lax}^{2}(Cat)$ (see the second question of Section~8) and
this is the main reason for us to stop at the notion of symmetric
bimonoidal intermuting categories.

This paper  gives, as a by-product, a complete formulation of a
fragment of linear derivations in classical and intuitionistic
propositional logic. Logic also helped the authors of \cite{DPInt}
to find that something is inappropriate in the unbalanced
treatment of units versus tensors in \cite{BFSV}. (Derivations of
sequents of the form $A\vee B\vdash A\kon B$ and $A\vee B\vdash
B\kon A$ are undesirable in logic.) We keep to the notation
$\vee,\kon,\bot,\top$ for tensors and units which is inspired by
logic. This is partly because at one point (see Section~4, Lemma
4.1) there is a reference to a coherence result from \cite{DP04}
where this notation is primary. Also, some easy lemmata in
Section~5 are taken over from \cite{DPInt}. Otherwise, this paper
is self-contained.

The categories we envisage are not just special n-fold monoidal
categories. It is not only the case that the difference would
appear in morphisms that involve the units, but the undesirable
morphisms mentioned above show that the axiomatization of $n$-fold
monoidal categories given in \cite{BFSV} leads to a
non-conservative extension of its fragment without the units. So,
the categories would be different in their unit-free fragments
too. Hence, to derive our coherence result, even for the unit-free
fragment, from the coherence result of \cite{BFSV} would be as to
derive Mac Lane's symmetric monoidal coherence of \cite{ML63} from
the cartesian coherence, which has a much simpler proof (see
\cite{K72a}, p.\ 129, where the result is announced, \cite{M80},
Theorem 2.2, \cite{TS96}, Theorem 8.2.3, p.\ 207, \cite{P02},
Section~7 and \cite{DP01}).

The first part of the paper is devoted to a coherence result for
symmetric bimonoidal intermuting categories. At one point, for
technical reasons, a strictification with respect to both
associativity and symmetry is used, and since the latter is not so
standard, although it is explained in details in \cite{DP04}, a
sketch of a proof why it actually works is given in Section~4. In
the second part of the paper, this coherence result is used to
establish that for every pair $(n,m)$ of natural numbers one can
iterate the reduced bar construction using first $n$ times the
monoidal structure given by $\vee$ and $\bot$, and then $m$ times
the monoidal structure given by $\kon$ and $\top$ of a symmetric
bimonoidal intermuting category $\cal C$ in order to obtain a lax
functor mapping  an $n+m$-tuple $(k_1,\ldots,k_{n+m})$ of natural
numbers, regerded as objects of the simplicial category, to ${\cal
C}^{k_1\cdot\ldots\cdot k_{n+m}}$.

The coherence result for symmetric bimonoidal intermuting
categories is already present in \cite{DPInt}. Although that paper
is not easy to read, this result, as well as the other coherence
results given there, is correct. The proof presented here is just
more self-contained and because of that, by our opinion, easier
for reading. However, the mathematical content remains the same.
So, the correct referring to this coherence result should go
through \cite{DPInt}.

Some parts of the paper may be skipped (Sections~4 and 6 are
optional) and for experts familiar with the work of Balteanu et
al.\ it is, perhaps, sufficient to see the definition of symmetric
bimonoidal intermuting categories (Section~2), then the statement
of a coherence result for these categories (Section~3, Theorem
3.1) and eventually Section~7, especially Lemma 7.2, which makes
this coherence result sufficient for the construction of a lax
functor with desired properties.

\section{Symmetric bimonoidal intermuting categories}

For $F$ a lax symmetric monoidal functor between ordinary
symmetric monoidal categories, let us call $F$ \emph{semistrong}
if the structural constraint with components
\[F(c)\otimes F(c')\str F(c\otimes c')
 \]
is an isomorphism and let us call $F$ \emph{normal} if the
constraint $I\str F(I)$ is an isomorphism.

\vspace{1ex}

\noindent{\sc Definition}\quad A \emph{symmetric bimonoidal
intermuting} category (shortly $SMI$ category) consists of the
following:

\vspace{.5ex}

1. a symmetric monoidal category $\langle{\cal M},\vee,\bot,
\check{\alpha},\check{\sigma},\check{\rho},\check{\lambda}\rangle$
(here $\check{\alpha}$, $\check{\sigma}$, $\check{\rho}$ and
$\check{\lambda}$ stay for associativity, symmetry, right and left
identity natural isomorphisms;
$\check{\alpha}^\str_{A,B,C}\!:A\vee(B\vee C)\str(A\vee B)\vee C$
has the inverse $\check{\alpha}^\rts_{A,B,C}$, etc.),

\vspace{.5ex}

2. a normal symmetric monoidal functor $\kon\!:{\cal M}\times
{\cal M}\str {\cal M}$,

\vspace{.5ex}

3. a semistrong symmetric monoidal functor $\top\!: 1\str {\cal
M}$,

\vspace{.5ex}

4. monoidal transformations $\hat{\alpha}$, $\hat{\sigma}$,
$\hat{\rho}$ and $\hat{\lambda}$ such that $\langle{\cal
M},\kon,\top,
\hat{\alpha},\hat{\sigma},\hat{\rho},\hat{\lambda}\rangle$
satisfies the pseudocommutative pseudomonoid conditions (i.e., the
equations of a symmetric monoidal category).

\vspace{1ex}

That $\kon$ is a lax symmetric monoidal functor means that there
is a natural transformation $\iota$ given by the family of arrows
\[
\iota_{A,B,C,D}\!:(A\kon B)\vee(C\kon D)\str(A\vee C)\kon(B\vee
D),
\]
and an arrow $\beta^\rts\!:\bot\str \bot\kon\bot$ such that the
following diagrams commute:

\vspace{3ex}

\begin{center}
\begin{picture}(200,80)
\put(0,10){\makebox(0,0){$(A\kon D)\vee((B\kon E)\vee(C\kon F))$}}
\put(200,10){\makebox(0,0){$((A\kon D)\vee(B\kon E))\vee(C\kon
F)$}}

\put(0,40){\makebox(0,0){$(A\kon D)\vee((B\vee C)\kon(E\vee F))$}}
\put(200,40){\makebox(0,0){$((A\vee B)\kon(D\vee E))\vee(C\kon
F)$}}

\put(0,70){\makebox(0,0){$(A\vee(B\vee C))\kon(D\vee(E\vee F))$}}
\put(200,70){\makebox(0,0){$((A\vee B)\vee C)\kon((D\vee E)\vee
F)$}}

\put(75,10){\vector(1,0){50}} \put(75,70){\vector(1,0){50}}
\put(0,18){\vector(0,1){14}} \put(200,18){\vector(0,1){14}}
\put(0,48){\vector(0,1){14}} \put(200,48){\vector(0,1){14}}

\put(100,7){\makebox(0,0)[t]{$\check{\alpha}^\str$}}
\put(100,73){\makebox(0,0)[b]{$\check{\alpha}^\str\kon\check{\alpha}^\str$}}
\put(-5,25){\makebox(0,0)[r]{$\mj\vee\iota$}}
\put(-5,55){\makebox(0,0)[r]{$\iota$}}
\put(205,25){\makebox(0,0)[l]{$\iota\vee\mj$}}
\put(205,55){\makebox(0,0)[l]{$\iota$}}

\put(100,40){\makebox(0,0){\small $(1)$}}

\end{picture}
\end{center}

\vspace{3ex}

\begin{center}
\begin{picture}(140,50)
\put(0,10){\makebox(0,0){$(A\kon C)\vee(B\kon D)$}}
\put(140,10){\makebox(0,0){$(B\kon D)\vee(A\kon C)$}}

\put(0,40){\makebox(0,0){$(A\vee B)\kon(C\vee D)$}}
\put(140,40){\makebox(0,0){$(B\vee A)\kon(D\vee C)$}}

\put(45,10){\vector(1,0){50}} \put(45,40){\vector(1,0){50}}
\put(0,18){\vector(0,1){14}} \put(140,18){\vector(0,1){14}}

\put(70,7){\makebox(0,0)[t]{$\check{\sigma}$}}
\put(70,43){\makebox(0,0)[b]{$\check{\sigma}\kon\check{\sigma}$}}
\put(-5,25){\makebox(0,0)[r]{$\iota$}}
\put(145,25){\makebox(0,0)[l]{$\iota$}}

\put(70,25){\makebox(0,0){\small $(2)$}}

\end{picture}
\end{center}

\vspace{3ex}

\begin{center}
\begin{picture}(130,50)
\put(0,10){\makebox(0,0){$(A\kon B)\vee\bot$}}
\put(130,10){\makebox(0,0){$A\kon B$}}

\put(0,40){\makebox(0,0){$(A\kon B)\vee(\bot\kon\bot)$}}
\put(130,40){\makebox(0,0){$(A\vee\bot)\kon(B\vee\bot)$}}

\put(35,10){\vector(1,0){70}} \put(45,40){\vector(1,0){40}}
\put(0,18){\vector(0,1){14}} \put(130,32){\vector(0,-1){14}}

\put(70,7){\makebox(0,0)[t]{$\check{\rho}^\str$}}
\put(65,44){\makebox(0,0)[b]{$\iota$}}
\put(-5,25){\makebox(0,0)[r]{$\mj\vee\beta^\rts$}}
\put(135,25){\makebox(0,0)[l]{$\check{\rho}^\str\kon\check{\rho}^\str$}}

\put(65,25){\makebox(0,0){\small $(3)$}}

\end{picture}
\end{center}

\noindent while the normality of $\kon$ means that there is an
arrow $\beta^\str\!:\bot\kon\bot\str\bot$ inverse to $\beta^\rts$.

That $\top$ is a lax symmetric monoidal functor means that there
are arrows $\tau^\str\!:\top\vee\top\str\top$ and
$\kappa\!:\bot\str\top$ such that the following diagrams commute:

\vspace{1ex}

\begin{center}
\begin{picture}(90,80)
\put(45,10){\makebox(0,0){$\top$}}

\put(0,40){\makebox(0,0){$\top\vee\top$}}
\put(90,40){\makebox(0,0){$\top\vee\top$}}

\put(0,70){\makebox(0,0){$\top\vee(\top\vee\top)$}}
\put(90,70){\makebox(0,0){$(\top\vee\top)\vee\top$}}

\put(0,65){\vector(0,-1){18}} \put(90,65){\vector(0,-1){18}}

\put(30,70){\vector(1,0){30}} \put(10,32){\vector(2,-1){30}}
\put(80,32){\vector(-2,-1){30}}

\put(45,73){\makebox(0,0)[b]{$\check{\alpha}^\str$}}
\put(-3,55){\makebox(0,0)[r]{$\mj\vee\tau^\str$}}
\put(93,55){\makebox(0,0)[l]{$\tau^\str\vee\mj$}}

\put(23,24){\makebox(0,0)[tr]{$\tau^\str$}}
\put(69,23){\makebox(0,0)[tl]{$\tau^\str$}}

\put(45,48){\makebox(0,0){\small $(4)$}}
\end{picture}
\end{center}

\vspace{1ex}

\begin{center}
\begin{picture}(90,50)
\put(45,10){\makebox(0,0){$\top$}}

\put(0,40){\makebox(0,0){$\top\vee\top$}}
\put(90,40){\makebox(0,0){$\top\vee\top$}}

\put(20,40){\vector(1,0){50}} \put(10,32){\vector(2,-1){30}}
\put(80,32){\vector(-2,-1){30}}

\put(45,43){\makebox(0,0)[b]{$\check{\sigma}$}}
\put(23,24){\makebox(0,0)[tr]{$\tau^\str$}}
\put(69,23){\makebox(0,0)[tl]{$\tau^\str$}}

\put(45,28){\makebox(0,0){\small $(5)$}}

\end{picture}
\end{center}

\vspace{1ex}

\begin{center}
\begin{picture}(90,50)
\put(45,10){\makebox(0,0){$\top$}}

\put(0,40){\makebox(0,0){$\top\vee\bot$}}
\put(90,40){\makebox(0,0){$\top\vee\top$}}

\put(20,40){\vector(1,0){50}} \put(10,32){\vector(2,-1){30}}
\put(80,32){\vector(-2,-1){30}}

\put(45,43){\makebox(0,0)[b]{$\mj\vee\kappa$}}
\put(23,24){\makebox(0,0)[tr]{$\check{\rho}^\str$}}
\put(69,23){\makebox(0,0)[tl]{$\tau^\str$}}

\put(45,28){\makebox(0,0){\small $(6)$}}

\end{picture}
\end{center}

\noindent while the semistrength of $\top$ means that there is an
arrow $\tau^\rts\!:\top\str \top\vee\top$ inverse to $\tau^\str$,
which enables us to shorten (4) into:

\vspace{1ex}

\begin{center}
\begin{picture}(90,50)
\put(45,10){\makebox(0,0){$\top\vee\top$}}

\put(0,40){\makebox(0,0){$\top\vee(\top\vee\top)$}}
\put(90,40){\makebox(0,0){$(\top\vee\top)\vee\top$}}

\put(30,40){\vector(1,0){30}} \put(10,32){\vector(2,-1){30}}
\put(80,32){\vector(-2,-1){30}}

\put(45,43){\makebox(0,0)[b]{$\check{\alpha}^\str$}}
\put(23,24){\makebox(0,0)[tr]{$\mj\vee\tau^\str$}}
\put(69,23){\makebox(0,0)[tl]{$\tau^\str\vee\mj$}}

\put(45,28){\makebox(0,0){\small $(4')$}}
\end{picture}
\end{center}

\vspace{1ex}

That $\hat{\alpha}$ is a monoidal transformation means that the
following diagrams commute:

\vspace{3ex}

\begin{center}
\begin{picture}(200,80)
\put(0,10){\makebox(0,0){$(A\vee D)\kon((B\vee E)\kon(C\vee F))$}}
\put(200,10){\makebox(0,0){$((A\vee D)\kon(B\vee E))\kon(C\vee
F)$}}

\put(0,40){\makebox(0,0){$(A\vee D)\kon((B\kon C)\vee(E\kon F))$}}
\put(200,40){\makebox(0,0){$((A\kon B)\vee(D\kon E))\kon(C\vee
F)$}}

\put(0,70){\makebox(0,0){$(A\kon(B\kon C))\vee(D\kon(E\kon F))$}}
\put(200,70){\makebox(0,0){$((A\kon B)\kon C)\vee((D\kon E)\kon
F)$}}

\put(75,10){\vector(1,0){50}} \put(75,70){\vector(1,0){50}}
\put(0,32){\vector(0,-1){14}} \put(200,32){\vector(0,-1){14}}
\put(0,62){\vector(0,-1){14}} \put(200,62){\vector(0,-1){14}}

\put(100,7){\makebox(0,0)[t]{$\hat{\alpha}^\str$}}
\put(100,73){\makebox(0,0)[b]{$\hat{\alpha}^\str\vee\hat{\alpha}^\str$}}
\put(-5,25){\makebox(0,0)[r]{$\mj\kon\iota$}}
\put(-5,55){\makebox(0,0)[r]{$\iota$}}
\put(205,25){\makebox(0,0)[l]{$\iota\kon\mj$}}
\put(205,55){\makebox(0,0)[l]{$\iota$}}

\put(100,40){\makebox(0,0){\small $(7)$}}

\end{picture}
\end{center}

\vspace{3ex}

\begin{center}
\begin{picture}(90,80)

\put(45,10){\makebox(0,0){$\bot$}}

\put(0,70){\makebox(0,0){$\bot\kon(\bot\kon\bot)$}}
\put(90,70){\makebox(0,0){$(\bot\kon\bot)\kon\bot$}}

\put(0,40){\makebox(0,0){$\bot\kon\bot$}}
\put(90,40){\makebox(0,0){$\bot\kon\bot$}}

\put(0,47){\vector(0,1){18}} \put(90,47){\vector(0,1){18}}

\put(30,70){\vector(1,0){30}} \put(40,17){\vector(-2,1){30}}
\put(50,17){\vector(2,1){30}}

\put(45,73){\makebox(0,0)[b]{$\hat{\alpha}^\str$}}
\put(-3,55){\makebox(0,0)[r]{$\mj\kon\beta^\rts$}}
\put(93,55){\makebox(0,0)[l]{$\beta^\rts\kon\mj$}}

\put(23,24){\makebox(0,0)[tr]{$\beta^\rts$}}
\put(69,23){\makebox(0,0)[tl]{$\beta^\rts$}}

\put(45,48){\makebox(0,0){\small $(8)$}}

\end{picture}
\end{center}
Since $\beta^\rts$ is invertible, (8) can be shorten into:

\vspace{1ex}

\begin{center}
\begin{picture}(90,50)

\put(45,10){\makebox(0,0){$\bot\kon\bot$}}

\put(0,40){\makebox(0,0){$\bot\kon(\bot\kon\bot)$}}
\put(90,40){\makebox(0,0){$(\bot\kon\bot)\kon\bot$}}

\put(30,40){\vector(1,0){30}} \put(40,17){\vector(-2,1){30}}
\put(50,17){\vector(2,1){30}}

\put(45,43){\makebox(0,0)[b]{$\hat{\alpha}^\str$}}
\put(23,24){\makebox(0,0)[tr]{$\mj\kon\beta^\rts$}}
\put(69,23){\makebox(0,0)[tl]{$\beta^\rts\kon\mj$}}

\put(45,28){\makebox(0,0){\small $(8')$}}

\end{picture}
\end{center}

\vspace{1ex}

That $\hat{\sigma}$ is a monoidal transformation means that the
following diagrams commute:

\vspace{3ex}

\begin{center}
\begin{picture}(140,50)
\put(0,10){\makebox(0,0){$(A\vee C)\kon(B\vee D)$}}
\put(140,10){\makebox(0,0){$(B\vee D)\kon(A\vee C)$}}

\put(0,40){\makebox(0,0){$(A\kon B)\vee(C\kon D)$}}
\put(140,40){\makebox(0,0){$(B\kon A)\vee(D\kon C)$}}

\put(45,10){\vector(1,0){50}} \put(45,40){\vector(1,0){50}}
\put(0,32){\vector(0,-1){14}} \put(140,32){\vector(0,-1){14}}

\put(70,7){\makebox(0,0)[t]{$\hat{\sigma}$}}
\put(70,43){\makebox(0,0)[b]{$\hat{\sigma}\vee\hat{\sigma}$}}
\put(-5,25){\makebox(0,0)[r]{$\iota$}}
\put(145,25){\makebox(0,0)[l]{$\iota$}}

\put(70,25){\makebox(0,0){\small $(9)$}}

\end{picture}
\end{center}

\vspace{3ex}

\begin{center}
\begin{picture}(90,50)

\put(45,10){\makebox(0,0){$\bot$}}

\put(0,40){\makebox(0,0){$\bot\kon\bot$}}
\put(90,40){\makebox(0,0){$\bot\kon\bot$}}

\put(20,40){\vector(1,0){50}} \put(40,17){\vector(-2,1){30}}
\put(50,17){\vector(2,1){30}}

\put(45,43){\makebox(0,0)[b]{$\hat{\sigma}$}}
\put(23,24){\makebox(0,0)[tr]{$\beta^\rts$}}
\put(69,23){\makebox(0,0)[tl]{$\beta^\rts$}}

\put(45,28){\makebox(0,0){\small $(10)$}}

\end{picture}
\end{center}

That $\hat{\rho}$ is a monoidal transformation means that the
following diagrams commute:

\vspace{3ex}

\begin{center}
\begin{picture}(130,50)
\put(0,10){\makebox(0,0){$A\vee B$}}
\put(130,10){\makebox(0,0){$(A\vee B)\kon\top$}}

\put(0,40){\makebox(0,0){$(A\kon\top)\vee(B\kon\top)$}}
\put(130,40){\makebox(0,0){$(A\vee B)\kon(\top\vee\top)$}}

\put(95,10){\vector(-1,0){70}} \put(45,40){\vector(1,0){40}}
\put(0,32){\vector(0,-1){14}} \put(130,32){\vector(0,-1){14}}

\put(70,7){\makebox(0,0)[t]{$\hat{\rho}^\str$}}
\put(65,44){\makebox(0,0)[b]{$\iota$}}
\put(-5,25){\makebox(0,0)[r]{$\hat{\rho}^\str\vee\hat{\rho}^\str$}}
\put(135,25){\makebox(0,0)[l]{$\mj\kon\tau^\str$}}

\put(65,25){\makebox(0,0){\small $(11)$}}

\end{picture}
\end{center}

\vspace{3ex}

\begin{center}
\begin{picture}(90,50)
\put(45,10){\makebox(0,0){$\bot$}}

\put(00,40){\makebox(0,0){$\bot\kon\bot$}}
\put(90,40){\makebox(0,0){$\bot\kon\top$}}

\put(20,40){\vector(1,0){50}} \put(40,17){\vector(-2,1){30}}
\put(80,32){\vector(-2,-1){30}}

\put(45,43){\makebox(0,0)[b]{$\mj\kon\kappa$}}
\put(23,24){\makebox(0,0)[tr]{$\beta^\rts$}}
\put(69,23){\makebox(0,0)[tl]{$\hat{\rho}^\str$}}

\put(45,28){\makebox(0,0){\small $(12)$}}

\end{picture}
\end{center}
(That $\hat{\lambda}$ is a monoidal transformation follows from
$\hat{\sigma}$ and $\hat{\rho}$ being so.)

Altogether, an $SMI$ category is a category $\cal M$ equipped with
two symmetric monoidal structures $\langle{\cal M},\vee,\bot,
\check{\alpha},\check{\sigma},\check{\rho},\check{\lambda}\rangle$
and $\langle{\cal M},\kon,\top,
\hat{\alpha},\hat{\sigma},\hat{\rho},\hat{\lambda}\rangle$, a
natural transformation $\iota$ given by the family of arrows
\[
\iota_{A,B,C,D}\!:(A\kon B)\vee(C\kon D)\str(A\vee C)\kon(B\vee
D),
\]
two isomorphisms $\beta^\rts\!:\bot\str\bot\kon\bot$,
$\tau^\str\!:\top\vee\top\str\top$, and an arrow
$\kappa\!:\bot\str\top$ such that the diagrams (1)-(12) commute.
(Note that the equations $(1)$ and $(7)$ are just unstrictified
versions of the external associativity condition and the internal
associativity condition of \cite{BFSV}.)

\prop{Lemma \thesection.1}{The equation
$\hat{\rho}^\str_\top\cirk(\kappa\kon\kappa)\cirk\beta^\rts=\kappa$
holds in every $SMI$ category.}

\dkz This follows from the commutative diagram

\begin{center}
\begin{picture}(240,50)
\put(80,10){\makebox(0,0){$\bot\kon\top$}}
\put(160,10){\makebox(0,0){$\bot$}}

\put(0,40){\makebox(0,0){$\bot$}}
\put(80,40){\makebox(0,0){$\bot\kon\bot$}}
\put(160,40){\makebox(0,0){$\top\kon\top$}}
\put(240,40){\makebox(0,0){$\top$}}

\put(10,30){\vector(3,-1){51}} \put(100,13){\vector(3,1){51}}
\put(180,13){\vector(3,1){51}} \put(80,32){\vector(0,-1){13}}
\put(10,40){\vector(1,0){51}} \put(100,40){\vector(1,0){41}}
\put(180,40){\vector(1,0){51}} \put(100,10){\vector(1,0){51}}

\put(33,20){\makebox(0,0)[t]{\scriptsize $\hat{\rho}^\rts_\bot$}}
\put(40,43){\makebox(0,0)[b]{\scriptsize $\beta^\rts$}}
\put(120,43){\makebox(0,0)[b]{\scriptsize $\kappa\kon\kappa$}}
\put(200,43){\makebox(0,0)[b]{\scriptsize $\hat{\rho}^\str_\top$}}
\put(120,7){\makebox(0,0)[t]{\scriptsize $\hat{\rho}^\str_\bot$}}
\put(77,25){\makebox(0,0)[r]{\scriptsize $\mj\kon\kappa$}}
\put(130,20){\makebox(0,0)[l]{\scriptsize $\kappa\kon\mj$}}
\put(210,20){\makebox(0,0)[l]{\scriptsize $\kappa$}}

\put(47,30){\makebox(0,0){\scriptsize $(12)$}}
\put(105,27){\makebox(0,0){\scriptsize $(bif.)$}}
\put(177,27){\makebox(0,0){\scriptsize $(nat.)$}}

\put(300,0){\makebox(0,0){$\dashv$}}

\end{picture}
\end{center}

\prop{Lemma \thesection.2}{The equation
$(\check{\rho}^\rts_\bot\kon
\check{\rho}^\rts_\bot)\cirk\beta^\rts=\iota\cirk(\beta^\rts \vee
\beta^\rts) \cirk \check{\rho}^\rts_\bot$ holds in every $SMI$
category.}

\dkz This follows from the commutative diagram

\begin{center}
\begin{picture}(240,50)
\put(0,10){\makebox(0,0){$\bot\vee\bot$}}
\put(80,10){\makebox(0,0){$(\bot\kon\bot)\vee\bot$}}
\put(200,10){\makebox(0,0){$(\bot\kon\bot)\vee(\bot\kon\bot)$}}

\put(0,40){\makebox(0,0){$\bot$}}
\put(80,40){\makebox(0,0){$\bot\kon\bot$}}
\put(200,40){\makebox(0,0){$(\bot\vee\bot)\kon(\bot\vee\bot)$}}

\put(0,32){\vector(0,-1){13}} \put(80,32){\vector(0,-1){13}}
\put(200,19){\vector(0,1){13}} \put(10,40){\vector(1,0){51}}
\put(100,40){\vector(1,0){51}} \put(20,10){\vector(1,0){30}}
\put(115,10){\vector(1,0){41}}

\put(40,43){\makebox(0,0)[b]{\scriptsize $\beta^\rts$}}
\put(125,43){\makebox(0,0)[b]{\scriptsize
$\check{\rho}^\rts_\bot\kon \check{\rho}^\rts_\bot$}}
\put(33,7){\makebox(0,0)[t]{\scriptsize $\beta^\rts\vee\mj$}}
\put(135,7){\makebox(0,0)[t]{\scriptsize $\mj\vee\beta^\rts$}}
\put(-3,25){\makebox(0,0)[r]{\scriptsize
$\check{\rho}^\rts_\bot$}}
\put(83,25){\makebox(0,0)[l]{\scriptsize
$\check{\rho}^\rts_{\bot\kon\bot}$}}
\put(203,25){\makebox(0,0)[l]{\scriptsize $\iota$}}

\put(140,25){\makebox(0,0){\scriptsize $(3)$}}
\put(40,25){\makebox(0,0){\scriptsize $(nat.)$}}

\put(300,0){\makebox(0,0){$\dashv$}}

\end{picture}
\end{center}

\prop{Proposition \thesection.3}{The following diagram commutes in
every $SMI$ category:}
\begin{center}
\begin{picture}(190,40)
\put(0,10){\makebox(0,0){$\bot\vee\bot$}}
\put(65,10){\makebox(0,0){$\bot$}}
\put(122,10){\makebox(0,0){$\top$}}
\put(190,10){\makebox(0,0){$\top\kon\top$}}

\put(0,40){\makebox(0,0){$(\top\kon\bot)\vee(\bot\kon\top)$}}
\put(190,40){\makebox(0,0){$(\top\vee\bot)\kon(\bot\vee\top)$}}

\put(22,10){\vector(1,0){30}} \put(78,10){\vector(1,0){30}}
\put(135,10){\vector(1,0){30}} \put(50,40){\vector(1,0){90}}
\put(0,32){\vector(0,-1){14}} \put(190,18){\vector(0,1){14}}

\put(37,7){\makebox(0,0)[t]{$\check{\rho}^\str$}}
\put(95,7){\makebox(0,0)[t]{$\kappa$}}
\put(150,7){\makebox(0,0)[t]{$\hat{\rho}^\rts$}}
\put(95,43){\makebox(0,0)[b]{$\iota$}}
\put(-5,25){\makebox(0,0)[r]{$\hat{\lambda}^\str\vee\hat{\rho}^\str$}}
\put(195,25){\makebox(0,0)[l]{$\check{\rho}^\rts\kon\check{\lambda}^\rts$}}

%\put(95,25){\makebox(0,0){\small $(13)$}}

\end{picture}
\end{center}

\dkz It suffices to show that the composite
\begin{center}
\begin{picture}(340,20)
\put(2,10){\makebox(0,0){$\bot$}}
\put(40,10){\makebox(0,0){$\bot\vee\bot$}}
\put(122,10){\makebox(0,0){$(\top\kon\bot)\vee(\bot\kon\top)$}}
\put(218,10){\makebox(0,0){$(\top\vee\bot)\kon(\bot\vee\top)$}}
\put(300,10){\makebox(0,0){$\top\kon\top$}}
\put(338,10){\makebox(0,0){$\top$}}

\put(9,10){\vector(1,0){17}} \put(55,10){\vector(1,0){24}}
\put(163,10){\vector(1,0){14}} \put(260,10){\vector(1,0){24}}
\put(315,10){\vector(1,0){17}}

\put(17,7){\makebox(0,0)[t]{\scriptsize $\check{\rho}^\rts$}}
\put(68,7){\makebox(0,0)[t]{\scriptsize
$\hat{\lambda}^\rts\vee\hat{\rho}^\rts$}}
\put(170,7){\makebox(0,0)[t]{\scriptsize $\iota$}}
\put(272,7){\makebox(0,0)[t]{\scriptsize
$\check{\rho}^\str\kon\check{\lambda}^\str$}}
\put(322,7){\makebox(0,0)[t]{\scriptsize $\hat{\rho}^\str$}}

\end{picture}
\end{center}
is equal to $\kappa$. This follows from the commutativity of the
diagram

\begin{center}
\begin{picture}(300,80)
\put(0,10){\makebox(0,0)[l]{$(\top\kon\bot)\vee(\bot\kon\top)$}}
\put(140,10){\makebox(0,0){$(\top\vee\bot)\kon(\bot\vee\top)$}}
\put(240,10){\makebox(0,0){$\top\kon\top$}}
\put(290,10){\makebox(0,0){$\top$}}

\put(0,40){\makebox(0,0)[l]{$\bot\vee\bot$}}
\put(100,40){\makebox(0,0){$(\bot\kon\bot)\vee(\bot\kon\bot)$}}
\put(200,40){\makebox(0,0){$(\bot\vee\bot)\kon(\bot\vee\bot)$}}
\put(290,40){\makebox(0,0){$\bot\kon\bot$}}

\put(10,70){\makebox(0,0)[l]{$\bot$}}
\put(200,70){\makebox(0,0){$\bot\kon\bot$}}

\put(82,10){\vector(1,0){16}} \put(183,10){\vector(1,0){40}}
\put(255,10){\vector(1,0){28}}

\put(13,32){\vector(0,-1){14}} \put(13,62){\vector(0,-1){14}}
\put(200,62){\vector(0,-1){14}}

\put(27,40){\vector(1,0){30}} \put(142,40){\vector(1,0){16}}
\put(242,40){\vector(1,0){34}}

\put(20,70){\vector(1,0){162}}

\put(220,70){\vector(3,-1){65}} \put(85,32){\vector(-2,-1){30}}
\put(185,32){\vector(-2,-1){30}} \put(280,32){\vector(-2,-1){30}}

\put(90,7){\scriptsize \makebox(0,0)[t]{$\iota$}}
\put(203,7){\scriptsize
\makebox(0,0)[t]{$\check{\rho}^\str\kon\check{\lambda}^\str$}}
\put(264,7){\scriptsize \makebox(0,0)[t]{$\hat{\rho}^\str$}}

\put(10,25){\scriptsize
\makebox(0,0)[r]{$\hat{\lambda}^\rts\vee\hat{\rho}^\rts$}}
\put(10,55){\scriptsize \makebox(0,0)[r]{$\check{\rho}^\rts$}}
\put(80,25){\scriptsize
\makebox(0,0)[l]{$(\kappa\!\kon\!\mj)\!\vee\!(\mj\!\kon\!\kappa)$}}
\put(180,25){\scriptsize
\makebox(0,0)[l]{$(\kappa\!\vee\!\mj)\!\kon\!(\mj\!\vee\!\kappa)$}}
\put(272,25){\scriptsize \makebox(0,0)[l]{$\kappa\kon\kappa$}}
\put(100,73){\scriptsize \makebox(0,0)[b]{$\beta^\rts$}}
\put(197,55){\scriptsize
\makebox(0,0)[r]{$\check{\rho}^\rts\kon\check{\rho}^\rts$}}
\put(250,65){\scriptsize \makebox(0,0)[l]{$\mj$}}

\put(42,43){\scriptsize
\makebox(0,0)[b]{$\beta^\rts\vee\beta^\rts$}}
\put(150,43){\scriptsize \makebox(0,0)[b]{$\iota$}}
\put(259,43){\scriptsize
\makebox(0,0)[b]{$\check{\rho}^\str\kon\check{\lambda}^\str$}}

\put(40,25){\scriptsize \makebox(0,0){$(12)$}}
\put(145,27){\scriptsize \makebox(0,0){$(nat.)$}}
\put(245,27){\scriptsize \makebox(0,0){$(nat.)$}}
\put(100,55){\scriptsize \makebox(0,0){(Lemma \thesection.2)}}

\end{picture}
\end{center}
where the composite along the top perimeter and going down is
$\kappa$ by Lemma \thesection.1. \qed

Our goal is to prove a coherence result for $SMI$ categories which
roughly says the following:

\vspace{2ex}

\noindent {\it Two canonical arrows $f,g\!:A\str B$ of an $SMI$
category are equal if either:}

\begin{itemize}

\item[] {\it the units $\bot$ and $\top$ do not ``essentially''
occur in $A$ and $B$, and $f$ and $g$ have the same graph} ({\it
defined analogously to the Kelly-Mac Lane graphs in}
\cite{KML71}), {\it or}

\item[] {\it $A$ and $B$ are isomorphic to $\bot$ or to $\top$.}

\end{itemize}
Since this result has to say something about the canonical
structure of an $SMI$ category, and this structure is equationally
presented, a precise formulation of our coherence result is given
in terms of an $SMI$ category freely generated by a set of
objects.

\section{Freely generated $SMI$ category}
Our category \cm\ (called ${\mathbf{SC^k}}_{\top,\bot}$ in
\cite{DPInt}), which is an $SMI$ category freely generated by an
infinite set $\cal P$ of propositional letters, is constructed as
follows:

The \emph{objects} of \cm\ are propositional formulae of the
language generated from $\cal P$, constants $\bot$ and $\top$,
with the binary connectives $\vee$ and $\kon$. The \emph{arrows}
of \cm\ are equivalence classes of \emph{arrow terms} generated
from \emph{primitive arrow terms} $\mj_A$,
$\check{\alpha}^\str_{A,B,C},\ldots$, $\kappa$, with the help of
$\cirk$, $\vee$ and $\kon$. These equivalence classes are taken
with respect to the smallest equivalence relation on arrow terms
which makes out of \cm\ an $SMI$ category. So, this equivalence
relation captures the equations of both symmetric monoidal
structures, naturality of $\iota$, isomorphism conditions for
$\beta$ and $\tau$, the equations brought by the commutative
diagrams (1)-(12) of the preceding section, and it is congruent
with respect to $\cirk$, $\vee$ and $\kon$.

Throughout this section we use the following terminology. We say
that an arrow term is an $\alpha$-\emph{term} if it is built from
identities and one occurrence of $\alpha$ with the help of $\vee$
and $\kon$ (see the definition of ``expanded instance of $a$''
given in \cite{ML63}). For example,
$\mj_A\kon(\hat{\alpha}^\str_{B,C,D}\vee\mj_E)$ is an
$\alpha$-term and we call $\hat{\alpha}^\str_{B,C,D}$ its
\emph{head}. We define analogously $\sigma$, $\rho$, $\lambda$,
$\iota$, $\beta$, $\tau$ and $\kappa$-terms and their heads. Note
that they are all composition free. We say that an arrow term
$f_n\cirk\ldots\cirk f_1\cirk\mj_A$ is a \emph{developed} arrow
term if each $f_i$ is $\alpha$, $\sigma$, $\rho$, $\lambda$,
$\iota$, $\beta$, $\tau$ or $\kappa$-term. It is easy to see that
every arrow term of \cm\ is equal to a developed one.

We say that an arrow term is \emph{defined by} $\alpha$ if it is
built from identities and $\alpha$'s (both $\check{\alpha}$'s and
$\hat{\alpha}$'s) with the help of $\cirk$, $\vee$ and $\kon$. We
say analogously that an arrow term is \emph{defined by} $\sigma$,
or by $\alpha$ and $\sigma$, etc.

Let us call $\nu$-terms all the $\rho$, $\lambda$, $\beta$ and
$\tau$-terms with superscripts $^\str$ in its heads. Let
$\rightharpoonup$ be a relation on objects of \cm\ defined by
$A\rightharpoonup B$ when there is a $\nu$-term $f\!:A\str B$.
Since $\rightharpoonup$ decreases the length of formulae, this
relation is noetherian, in the terminology of \cite{H80}. We can
also prove the following lemma which is analogous to a lemma that
implicitly occurs in the proof of monoidal coherence given in
\cite{ML63} (Section~3) or in \cite{ML71} (Section VII.2).

\prop{Lemma \thesection.1}{The relation $\rightharpoonup$ is
locally confluent and this is justified by commutative diagrams of
$\nu$-terms.}

\dkz In all the possible cases for a pair of $\nu$-terms
$(f\!:A\str B,g\!:A\str C)$ that are not equal in \cm, we use
either that $\vee$ and $\kon$ are bifunctors, or that $\rho$'s and
$\lambda$'s are natural in order to find a pair of $\nu$-terms
$(g'\!:B\str D,f'\!:C\str D)$ such that $g'\cirk f=f'\cirk g$.
\qed

Since $\rightharpoonup$ is noetherian and locally confluent every
object $A$ of \cm\ has a unique normal form which we denote by
$\nu(A)$, and we say that $A$ \emph{reduces by} $\nu$ to $\nu(A)$.
If no letter occurs in $A$, then $\nu(A)$ is either $\bot$ or
$\top$.

We call the arrow terms defined by $\rho$, $\lambda$, $\beta$ and
$\tau$, $\mN_{\top,\bot}$-\emph{terms} as in \cite{DPInt}, and
when all the superscripts are $^\str$ we call them
\emph{directed}. Since every directed $\mN_{\top,\bot}$-term is
equal to a developed one (i.e.\ to a composition of an identity
and some $\nu$-terms) as a corollary of Lemma \thesection.1 we
have the following.

\prop{Lemma \thesection.2}{If $f,g\!:A\str\nu(A)$ are two directed
$\mN_{\top,\bot}$-terms, then $f=g$.}

Next we can prove the following.

\prop{Lemma \thesection.3}{Every diagram of
$\mN_{\top,\bot}$-terms is commutative.}

\dkz This is established in the same way as the coherence result
for monoidal categories in \cite{ML63} or in \cite{ML71} (Section
VII.2) by relying on Lemma \thesection.2.

\prop{Lemma \thesection.4}{If no letter occurs in $A$ and $B$ then
for every pair $f,g\!:A\str B$ of arrow terms defined by $\rho$,
$\lambda$, $\beta$, $\tau$ and $\kappa$, we have $f=g$.}

\dkz We establish first that every arrow term defined by $\rho$,
$\lambda$, $\beta$, $\tau$ and $\kappa$ is either equal to an
$\mN_{\top,\bot}$-term or it is equal to a term of the form
$h''\cirk\kappa\cirk h'$ for $\mN_{\top,\bot}$-terms $h'$ and
$h''$. To do this, we rely on the equations (6), (12), the
following naturality conditions
\begin{tabbing}
\hspace{5em}\=$\kappa\vee\mj_\bot=\check{\rho}^\rts_\top\cirk
\kappa\cirk\check{\rho}^\str_\bot$,\hspace{3em}\=$\mj_\bot\vee\kappa
=\check{\lambda}^\rts_\top\cirk
\kappa\cirk\check{\lambda}^\str_\bot$,
\\[2ex]
\>$\kappa\kon\mj_\top=\hat{\rho}^\rts_\top\cirk
\kappa\cirk\hat{\rho}^\str_\bot$,\>$\mj_\top\kon\kappa
=\hat{\lambda}^\rts_\top\cirk \kappa\cirk\hat{\lambda}^\str_\bot$,
\end{tabbing}
and the fact that there are no arrow terms of the form
$\kappa\cirk h\cirk\kappa$. If $f$ is equal to an
$\mN_{\top,\bot}$-term then $A$ is isomorphic to $B$ and so $g$
must be equal to an $\mN_{\top,\bot}$-term too, and vice versa. It
only remains to apply Lemma \thesection.3. \qed

Let $i_A\!:A\str\nu(A)$ be a directed $\mN_{\top,\bot}$-term. By
Lemma \thesection.2 we know that any choice of the term $i_A$
leads to the same isomorphism of \cm. We can prove the following.

\prop{Lemma \thesection.5}{If no letter occurs in $A$ and $B$ then
every arrow $f\!:A\str B$ may be defined by $\rho$, $\lambda$,
$\beta$, $\tau$ and $\kappa$.}

\dkz We rely on the equations of symmetric monoidal categories,
the equations (3), (4'), (5), (8'), (10), (11), Proposition 2.3,
the naturality conditions, and the fact that the arrow terms $i_A$
and $i^{-1}_A$ are $\mN_{\top,\bot}$-terms, to eliminate the
presence of $\alpha$'s, $\sigma$'s and $\iota$'s. For example, we
have

\vspace{1ex}

\begin{center}
\begin{picture}(180,90)
\put(180,80){\makebox(0,0){$(\bot\kon(\top\kon(\bot\kon\bot)))\kon\bot$}}
\put(180,40){\makebox(0,0){$(\bot\kon\bot)\kon\bot$}}

\put(90,10){\makebox(0,0){$\bot\kon\bot$}}

\put(0,80){\makebox(0,0){$\bot\kon((\top\kon(\bot\kon\bot))\kon\bot)$}}
\put(0,40){\makebox(0,0){$\bot\kon(\bot\kon\bot)$}}

\put(35,40){\vector(1,0){110}} \put(60,80){\vector(1,0){60}}
\put(0,72){\vector(0,-1){24}} \put(180,48){\vector(0,1){24}}
\put(10,30){\vector(4,-1){60}} \put(110,15){\vector(4,1){60}}

\put(-5,60){\makebox(0,0)[r]{\small $\mj\kon(i\kon\mj)$}}
\put(185,60){\makebox(0,0)[l]{\small $(\mj\kon i^{-1})\kon\mj$}}
\put(90,83){\makebox(0,0)[b]{\small $\hat{\alpha}^\str$}}
\put(90,43){\makebox(0,0)[b]{\small $\hat{\alpha}^\str$}}
\put(30,20){\makebox(0,0)[tr]{\small $\mj\kon\beta^\str$}}
\put(150,20){\makebox(0,0)[tl]{\small $\beta^\rts\kon\mj$}}

\put(90,63){\makebox(0,0){\small (\emph{nat})}}
\put(90,25){\makebox(0,0){\small $(8')$}}

\put(250,0){\makebox(0,0){$\dashv$}}

\end{picture}
\end{center}

\vspace{1ex}

As a direct consequence of Lemmata \thesection.4 and \thesection.5
we have:

\prop{Lemma \thesection.6}{If no letter occurs in $A$ and $B$ then
for every $f,g\!:A\str B$ we have $f=g$.}

Here is the explanation what we meant by not ``essential''
occurrence of the units in an object. We say that an object $A$ of
\cm\ is $\bot$-\emph{pure} when there is no occurrence of $\bot$
in $\nu(A)$. It is easy to see that $A$ is not $\bot$-pure iff
either $\nu(A)=\bot$ or there is a conjunction in $A$ (by a
conjunction in $A$ we mean a subformula of the form $B\kon C$)
such that one of its conjuncts reduces by $\nu$ to $\bot$ and a
letter occurs in the other. We define analogously a
$\top$-\emph{pure} object of \cm\ and derive an analogous
characterization. An object of \cm\ is \emph{pure} when it is both
$\bot$-pure and $\top$-pure. Pure objects play here a role similar
to the role of proper shapes in the symmetric monoidal closed
coherence proved by Kelly and Mac Lane in \cite{KML71}.

\prop{Lemma \thesection.7}{Let $f\!:A\str B$ be an arrow of \cm.
If $A$ is $\bot$-pure, then $B$ is $\bot$-pure, and if $B$ is
$\top$-pure, then $A$ is $\top$-pure.}

\dkz Since $f$ may be represented by a developed term it is
sufficient to verify the lemma for $\alpha$, $\sigma$, $\rho$,
$\lambda$, $\iota$, $\beta$, $\tau$ and $\kappa$-terms. The only
interesting case is when $f$ is a $\iota$-term.

Suppose $B$ is not $\bot$-pure. By using the above-mentioned
characterization of such objects of \cm, we have two
possibilities. If $\nu(B)=\bot$ then we easily conclude that
$\nu(A)=\bot$ too. If there is a conjunction in $B$ such that one
of its conjuncts is reduced by $\nu$ to $\bot$ and a letter occurs
in the other conjunct, then we obviously have the same situation
in $A$, except in the  case when this conjunction is the target of
the head $\iota_{E,F,G,H}\!:(E\kon F)\vee(G\kon H)\str(E\vee
G)\kon(F\vee H)$ of $f$. If $\nu(E\vee G)$ is $\bot$ and there is
a letter in $F\vee H$, then $\nu(E)=\nu(G)=\bot$ and there is a
letter in either $F$ or $G$. So, $A$ is not $\bot$-pure. This is
sufficient for the first implication and we proceed analogously
for the second implication of the lemma. \qed

\prop{Corollary}{If $f\!:A\str B$ and $g\!:B\str C$ are arrows of
\cm\ such that $A$ and $C$ are pure, then $B$ is pure.}

\vspace{-2ex}

\prop{Lemma \thesection.8}{If $f\!:A\str B$ is an arrow term such
that $A$ and $B$ are pure, then there is an arrow term
$f'\!:\nu(A)\str\nu(B)$ such that $\rho$, $\lambda$, $\beta$,
$\tau$, $\kappa$, $\top$ and $\bot$ do not occur in $f'$ and
\[f=i^{-1}_B\cirk f'\cirk i_A.\]}

\vspace{-2ex}

\dkz By the corollary of Lemma \thesection.7, it is sufficient to
verify the lemma for $\alpha$, $\sigma$, $\rho$, $\lambda$,
$\iota$, $\beta$, $\tau$ and $\kappa$-terms. If $f$ is an
$\mN_{\top,\bot}$-term then $\nu(A)=\nu(B)$ and by Lemma
\thesection.3 we have $f=i^{-1}_B\cirk i_A$.

If $f$ is an $\alpha$-term whose head is
$\hat{\alpha}^\str_{C,D,E}$, then by the following naturality
diagram

\begin{center}
\begin{picture}(130,50)
\put(0,10){\makebox(0,0){$\nu(X)\kon(\nu(Y)\kon\nu(Z))$}}
\put(130,10){\makebox(0,0){$(\nu(X)\kon\nu(Y))\kon\nu(Z)$}}

\put(0,40){\makebox(0,0){$X\kon(Y\kon Z)$}}
\put(130,40){\makebox(0,0){$(X\kon Y)\kon Z$}}

\put(50,10){\vector(1,0){30}} \put(35,40){\vector(1,0){60}}
\put(0,32){\vector(0,-1){14}} \put(130,18){\vector(0,1){14}}

\put(67,7){\makebox(0,0)[t]{\small $\hat{\alpha}^\str$}}
\put(67,44){\makebox(0,0)[b]{\small $\hat{\alpha}^\str$}}
\put(-5,25){\makebox(0,0)[r]{\small $i\kon(i\kon i)$}}
\put(135,25){\small \makebox(0,0)[l]{$(i^{-1}\kon i^{-1})\kon
i^{-1}$}}

\put(65,25){\makebox(0,0){\small (\emph{nat})}}

\end{picture}
\end{center}
we may assume that the indices $C$, $D$ and $E$ are already
reduced by $\nu$. By the assumption that $A$ and $B$ are pure we
have the following cases:

(1) the units do not occur in $C$, $D$ and $E$; hence we are
already done,

(2) one of $C$, $D$ or $E$ is $\top$; we are done by the following
commutative diagram delivered by the second monoidal structure
(here we assume $C=\top$ and we proceed analogously when $D=\top$
or $E=\top$),

\begin{center}
\begin{picture}(130,50)
\put(65,10){\makebox(0,0){$D\kon E$}}

\put(0,40){\makebox(0,0){$\top\kon(D\kon E)$}}
\put(130,40){\makebox(0,0){$(\top\kon D)\kon E$}}

\put(35,40){\vector(1,0){60}} \put(10,32){\vector(2,-1){30}}
\put(90,17){\vector(2,1){30}}

\put(65,43){\makebox(0,0)[b]{\small $\hat{\alpha}^\str$}}
\put(23,24){\makebox(0,0)[tr]{\small $\hat{\lambda}^\str$}}
\put(109,23){\makebox(0,0)[tl]{\small
$\hat{\lambda}^\rts\kon\mj$}}

\end{picture}
\end{center}

(3) $C=D=E=\bot$; we are done by relying on the equation~(8').

\vspace{1ex}

\noindent The situation is quite similar with the other $\alpha$
and $\sigma$-terms.

If $f$ is a $\iota$-term then again by naturality we may assume
that all the indices of the head of $f$ are reduced by $\nu$. It
is not possible that only one of its indices is reduced to $\bot$
or to $\top$ since then $A$ or $B$ is not pure. If two of its
indices are $\bot$ or $\top$ while the units do not occur in the
remaining two indices, then by the assumption that $A$ and $B$ are
pure, we may eliminate this $\iota$ by applying the equations (3)
or (11). Situation is analogous when three indices of $\iota$ are
$\bot$ or $\top$ and the forth is not. If all the indices of
$\iota$ are $\bot$ or $\top$ then we have two cases: either we
apply the equations (3) or (11) to eliminate $\iota$, or we apply
Proposition 2.3 to reduce $\iota$ to $\kappa$ and we deal with the
new occurrence of $\kappa$ as in the following last case for $f$.

If $f$ is a $\kappa$-term. Since $A$ and $B$ are pure, $f$ is not
just $\kappa$, so the head of $f$ is in the immediate scope of
$\vee$ or $\kon$. If $\mj_E\vee\kappa$ is a subterm of $f$, then
since $A$ and $B$ are pure, no letter occurs in $E$ and again we
may assume that $E$ is already reduced by $\nu$ to $\bot$ or
$\top$. If $E$ is $\top$ then we use the equation (6) to eliminate
$\kappa$. If $E$ is $\bot$ then we apply the naturality equation
$\mj_\bot\vee\kappa =\check{\lambda}^\rts_\top\cirk
\kappa\cirk\check{\lambda}^\str_\bot$, mentioned in the proof of
Lemma \thesection.4. This equation does not eliminate $\kappa$ but
it replaces a $\kappa$-term of a greater complexity by a
$\kappa$-term of lower complexity and by induction $\kappa$ will
be eliminated.

We proceed analogously in all the other possible cases for a
$\kappa$ term $f$ relying on equations (6), (12) or the remaining
naturality conditions mentioned in the proof of Lemma
\thesection.4. \qed

Formulations of coherence results that are not of the form ``all
diagrams commute'' usually require a notion of graph or diagram
associated to every canonical arrow of the structure for which the
result is formulated. Such a coherence result says that ``if
$f,g:\!A\str B$ have the same graph then $f=g$''. Sometimes these
graphs correspond to relations, functions, bijections, or like in
the case of Kelly-Mac Lane graphs, to Brauerian diagrams (see
\cite{DP07}, Section 2.3 and references therein). They are closely
related to the notion of generality formalized by Lambek in
\cite{L68} and \cite{L69} (see also \cite{DP03a}).

We say that an object of \cm\ is \emph{diversified} if every
letter occurs in it at most once. By induction on the complexity
of arrow term it can be shown that it is an instance of an arrow
term whose source and target are diversified and share the
letters. The graph associated to an arrow term $f$ corresponds to
the bijection between the letters in the source and target of a
``diversified'' arrow term whose instance is $f$. So, for \cm, we
can conclude that: ``if $f,g:\!A\str B$ have the same graph, then
$f=g$'' is equivalent to ``if $A$ and $B$ are diversified, then
there is at most one arrow $f\!:A\str B$''. From left to right
this is trivial and for the other direction, we use $f'$ and $g'$
with the same diversified source and target whose instances are
$f$ and $g$ respectively. By the assumption we obtain $f'=g'$ and
hence $f=g$ (just in the proof of $f'=g'$ use the same
substitution of letters needed to obtain $f$ from $f'$ and $g$
from $g'$).

So, our coherence for $SMI$ categories (called Restricted
Symmetric Bimonoidal Intermuting Coherence in \cite{DPInt}) is
formulated as follows:

\prop{Theorem \thesection.1}{If $A$ and $B$ are either pure and
diversified or no letter occurs in them, then there is at most one
arrow $f\!:A\str B$ in \cm.}

\noindent One part of the theorem is established by Lemma
\thesection.6. By Lemma \thesection.8 we have reduced the rest of
the theorem to the case when the units do not occur in $A$ and
$B$, and $f$ and $g$ are defined by $\alpha$, $\sigma$ and
$\iota$. So, to complete the proof of Theorem \thesection.1 it is
sufficient to prove a coherence result for categories like $SMI$
categories but without units, which we call as in \cite{DPInt},
\emph{symmetric biassociative intermuting} ($SAI$) categories. The
canonical structure of $SAI$ categories is given by two
biendofunctors $\vee$ and $\kon$, natural isomorphisms given by
associativities $\alpha$ and symmetries $\sigma$ that satisfy Mac
Lane's pentagonal and hexagonal conditions, and a natural
transformation $\iota$ satisfying the coherence conditions given
by the diagrams (1), (2), (7) and (9).

This coherence result is formulated in terms of the category that
should be called $SAI(\cal P)$, but we call it here simply \ca\
because it is an auxiliary category and has two modifications,
namely $\ca'$ and $\ca^{st}$, which we use later for our proof.
The category \ca\ is freely generated $SAI$ category by the same
set $\cal P$ of generators as \cm. The construction of \ca\ is
analogous to the construction of \cm\ given at the beginning of
this section. So, our auxiliary coherence result is the following:

\prop{Theorem \thesection.2}{If  $A$ and $B$ are diversified, then
there is at most one arrow $f\!:A\str B$ in \ca.}

\noindent The following two sections contain a proof of this
theorem.

\section{A note on strictification}
In order to provide an easier record of equations of arrow terms
in the proof of Theorem 3.2 we will replace our category \ca\ by a
symmetric biassociative  intermuting category in which
associativity and symmetry arrows are identities. Strictification
under associativity is a standard procedure in coherence results.
For example, this is how Mac Lane reduced his proof of symmetric
monoidal coherence in \cite{ML63} to the standard presentation of
symmetric groups by generators and relations. However,
strictification under symmetry is not so standard and it may cause
a suspicion. (A reference where it is used implicitly is
\cite{E66}.) Although various strictifications, including this
with respect to symmetry, are thoroughly investigated in
\cite{DP04}, Chapter~3 and \S\S4.7, 7.6-8, 8.4, we briefly pass
through such a strictification of our category \ca.

Note first that if we factor the arrow terms of \ca\ by the new
equations
\[
\check{\sigma}_{A,A}=\mj_{A\vee A},\quad
\hat{\sigma}_{A,A}=\mj_{A\kon A}
\]
obtaining a new category $\ca'$ with the same objects as \ca, the
full subcategories of \ca\ and $\ca'$ on diversified objects are
the same. This is because we can easily establish that for every
pair of arrow terms $f,g\!:A\str B$, if $f=g$ in $\ca'$ and $f\neq
g$ in \ca, then $A$ and $B$ are not diversified. Since the objects
$A$ and $B$ are diversified in Theorem 3.2, we can replace the
category \ca\ in the formulation of that theorem by the category
$\ca'$ without losing its strength. We use this fact later on.

Let the arrow terms defined by associativities $\alpha$ and
symmetries $\sigma$ (cf.\ the beginning of the preceding section)
be called $S$-\emph{terms}. Then we have the following result from
\cite{DP04}, \S6.5.

\prop{Lemma \thesection.1}{Every diagram of $S$-terms commutes in
$\ca'$.}

\noindent This fact together with the property that every $S$-term
represents an isomorphism of $\ca'$ is sufficient for our
strictification of $\ca'$ with respect to its associative and
symmetric structures. Roughly speaking, we can further factor the
arrow terms so that associativity and symmetry natural
transformations become identity natural transformations. Of,
course, this makes some identifications among the objects of
$\ca'$ too.

We define a relation $\equiv$ on the set of objects of $\ca'$
(which are the same as the objects of \ca) in the following way.
Let $A\equiv B$ iff there is an $S$-term $f\!:A\str B$. Since
$\mj_A$ is an $S$-term, every $S$-term represents an isomorphism
whose inverse may be represented by an $S$-term, and the
composition of two $S$-terms is an $S$-term, we have that $\equiv$
is an equivalence relation. Let $\lz A\dz$ denotes the equivalence
class with respect to $\equiv$ of an object $A$ of $\ca'$.

Since the objects of \ca, and hence of $\ca'$, are propositional
formulae of the language generated from $\cal P$, with the binary
connectives $\vee$ and $\kon$, they correspond to planar binary
trees with elements of $\cal P$ in the leaves and $\vee$ or $\kon$
in the vertices (see \cite{MT10}, Section 2.1, for the definition
of planar tree). Our relation $\equiv$ is such that if $A$
corresponds to a planar tree $T$, then $\lz A\dz$ corresponds to
the non-planar tree obtained from $T$ by omitting the associated
linear ordering, and by contraction of every edge having the same
connective in its ends (see \cite{BM76}, Section 2.4, for the
definition of the operation of contraction of an edge). Hence, we
can denote (not in a unique way) the equivalence class $\lz A\dz$
by deleting from the formula $A$ parenthesis tied to $\vee$ in the
immediate scope of another $\vee$ and the same for $\kon$. For
example, the equivalence class $\lz(p\kon q)\kon((p\vee r)\vee
p)\dz$ is denoted by $p\kon q\kon(p\vee r\vee p)$, and the same
equivalence class may be denoted by $q\kon(r\vee p\vee p)\kon p$
or by $(p\vee p\vee r)\kon p\kon q$, etc.

We call $\lz A\dz$ a \emph{form multiset} (see \cite{DP04},
\S7.7), in particular, when $A$ is diversified we call $\lz A\dz$
a \emph{form set}. We use $S$, $T$, $U$, $V$, $W$, $X$, $Y$ and
$Z$, possible with indices, for form multisets and form sets.

Note that if $A_1\equiv A_2$ and $B_1\equiv B_2$, then $A_1\vee
B_1\equiv A_2\vee B_2$ and $A_1\kon B_1\equiv A_2\kon B_2$, hence
we may define the operations $\vee$ and $\kon$ on form multisets
as
\[
\lz A\dz\vee\lz B\dz=_{df}\lz A\vee B\dz,\quad \lz A\dz\kon\lz
B\dz=_{df}\lz A\kon B\dz.
\]

Let $\ca^{st}$ be a category built out of syntactical material,
starting from the same set $\cal P$ of generators as in the case
of \cm, \ca\ and $\ca'$, whose objects are the form multisets. The
only \emph{primitive arrow terms} of $\ca^{st}$ are of the form
\[
\mj_S\!:S\str S,\quad\mbox{or}\quad\iota_{[S,T,U,V]}\!:(S\kon
T)\vee(U\kon V)\str(S\vee U)\kon(T\vee V),
\]
where $[S,T,U,V]$ is an abbreviation for the set
\[
\{(S,T,U,V),(T,S,V,U),(U,V,S,T),(V,U,T,S)\}.
\]

Hence, $\iota_{[S,T,U,V]}$, $\iota_{[T,S,V,U]}$
$\iota_{[U,V,S,T]}$, $\iota_{[V,U,T,S]}$ are the same primitive
arrow term which prevents us for having many primitive arrow terms
representing the same arrow of $\ca^{st}$. Moreover, the
strictified versions of the equations (2) and (9) are now
incorporated in our notation, and when we draw the arrow
$\iota\!:(S\kon T)\vee(U\kon V)\str(S\vee U)\kon(T\vee V)$ in a
diagram, one can form the index of $\iota$ in a unique way.

The \emph{arrows} of $\ca^{st}$ are equivalence classes of
\emph{arrow terms} generated from primitive arrow terms with the
help of $\cirk$, $\vee$ and $\kon$. These equivalence classes are
taken with respect to the smallest equivalence relation on arrow
terms which makes out of $\ca^{st}$ a strict associative and
strict symmetric $SAI$ category. So, this equivalence relation is
congruent with respect to $\cirk$, $\vee$ and $\kon$, and it
captures the assumptions that $\vee$ and $\kon$ are
biendofunctors, the following equations
\begin{tabbing}
\hspace{5em}\=$s\vee(t\vee u)=(s\vee t)\vee
u$,\hspace{5em}\=$s\kon(t\kon u)=(s\kon t)\kon u$,
\\[1ex]
\>$s\vee t=t\vee s$,\>$s\kon t=t\kon s$,
\end{tabbing}
(which are the rudiments of naturality conditions for
associativity and symmetry), naturality of $\iota$, and the
equations brought by the following commutative diagrams:

\vspace{1ex}

\begin{center}
\begin{picture}(200,80)
\put(100,10){\makebox(0,0){$(U\kon X)\vee(V\kon Y)\vee(W\kon Z)$}}

\put(0,40){\makebox(0,0){$(U\kon X)\vee((V\vee W)\kon(Y\vee Z))$}}
\put(200,40){\makebox(0,0){$((U\vee V)\kon(X\vee Y))\vee(W\kon
Z)$}}

\put(100,70){\makebox(0,0){$(U\vee V\vee W)\kon(X\vee Y\vee Z)$}}

\put(75,19){\vector(-4,1){50}} \put(125,19){\vector(4,1){50}}
\put(25,49){\vector(4,1){50}} \put(175,49){\vector(-4,1){50}}

\put(35,23){\makebox(0,0)[r]{\small $\mj\vee\iota$}}
\put(35,59){\makebox(0,0)[r]{\small $\iota$}}
\put(165,23){\makebox(0,0)[l]{\small $\iota\vee\mj$}}
\put(165,59){\makebox(0,0)[l]{\small $\iota$}}

\put(100,40){\makebox(0,0){\small $(1s)$}}

\end{picture}
\end{center}

\vspace{3ex}

\begin{center}
\begin{picture}(200,80)
\put(100,10){\makebox(0,0){$(U\vee X)\kon(V\vee Y)\kon(W\vee Z)$}}

\put(0,40){\makebox(0,0){$(U\vee X)\kon((V\kon W)\vee(Y\kon Z))$}}
\put(200,40){\makebox(0,0){$((U\kon V)\vee(X\kon Y))\kon(W\vee
Z)$}}

\put(100,70){\makebox(0,0){$(U\kon V\kon W)\vee(X\kon Y\kon Z)$}}

\put(25,32){\vector(4,-1){50}} \put(175,32){\vector(-4,-1){50}}
\put(75,62){\vector(-4,-1){50}} \put(125,62){\vector(4,-1){50}}

\put(35,23){\makebox(0,0)[r]{$\mj\kon\iota$}}
\put(35,59){\makebox(0,0)[r]{$\iota$}}
\put(165,23){\makebox(0,0)[l]{$\iota\kon\mj$}}
\put(165,59){\makebox(0,0)[l]{$\iota$}}

\put(100,40){\makebox(0,0){\small $(7s)$}}

\end{picture}
\end{center}
This concludes the definition of $\ca^{st}$.

The categories $\ca'$ and $\ca^{st}$ are equivalent via functors
that preserve the $SAI$ structure. Here is just a sketch of the
proof. We define two functors $H_{\cal G}\!:\ca'\str\ca^{st}$ and
$H\!:\ca^{st}\str\ca'$ in the following way. Let $H_{\cal
G}A=_{df}\lz A\dz$, and let $H_{\cal G}f$ be obtained from the
arrow term $f$ by replacing every $S$-term in it by $\mj$ indexed
by the equivalence class of the source and the target of this
term, and by replacing every $\iota_{A,B,C,D}$ in it by
$\iota_{[\lz A\dz,\lz B\dz,\lz C\dz,\lz D\dz\,]}$. It is not
difficult to verify that $H_{\cal G}$ is indeed a functor, i.e.\
that if $f=g$ in $\ca'$ then $H_{\cal G}f=H_{\cal G}g$ in
$\ca^{st}$.

On the other hand, to define $H\!:\ca^{st}\str\ca'$ we have first
to choose a formula $A_H$ in each equivalence class $\lz A\dz$. By
Lemma \thesection.1, there is a unique arrow $\varphi_A\!:A_H\str
A$ of $\ca'$ represented by an $S$-term. We define

\begin{tabbing}
\centerline{$H\lz A\dz=_{df}A_H$,}
\\[2ex]
\hspace{7em}\=$H\mj_S=_{df}\mj_{HS},$\hspace{5em}\=$H(t\cirk
s)=_{df}Ht\cirk Hs$,
\\[2ex]
\hspace{1em}$H\iota_{[S,T,U,V]}=_{df}\varphi^{-1}_{(HS\vee
HU)\kon(HT\vee HV)}\cirk\iota_{HS,HT,HU,HV}\cirk\varphi_{(HS\kon
HT)\vee(HU\kon HV)}$,
\\[2ex]
\hspace{1em}$H(s\vee t)=\varphi^{-1}_{HS_2\vee HT_2}\cirk(Hs\vee
Ht)\cirk\varphi_{HS_1\vee HT_1}$,\hspace{1em} for $s\!:S_1\str
S_2$, $t\!:T_1\str T_2$,
\\*[1ex] and the same for $\vee$ replaced by $\kon$.
\end{tabbing}

It can be easily checked that this definition is correct and that
so defined $H$ is indeed a functor. It is straightforward that
$H_{\cal G}\cirk H$ is the identity functor on $\ca^{st}$ and one
can verify that $\varphi$, defined as above, is a natural
isomorphism from $H\cirk H_{\cal G}$ to the identity functor on
$\ca'$. (Details of the proof, but in more general context, are
given in \cite{DP04}, \S 3.2.) Hence, $\ca'$ and $\ca^{st}$ are
equivalent via $H_{\cal G}$ and $H$. Following the terminology of
\cite{DP04}, functor $H_{\cal G}$ strictly preserves $SAI$
structure and $H$ is just strong with respect to this structure.

As a consequence of this equivalence and the fact that \ca\ and
$\ca'$ have the same full subcategories on diversified objects, we
have that the following coherence result is sufficient for Theorem
3.2.

\prop{Proposition \thesection.2}{If $X$ and $Y$ are form sets,
then there is at most one arrow $t\!:X\str Y$ in $\ca^{st}$.}

As we said at the beginning of this section, the strictification
of \ca\ enables us to record our derivations in the proof of this
proposition, and there are no other reasons, except these
technical, for this step. Note that one can always decorate the
arrow terms of $\ca^{st}$ (using the functor $H$) by lengthy
compositions of $S$-terms to get back into a rather natural
environment given by the category~\ca.

\section{Proof of Proposition~4.2}
In this section we are interested only in form sets (i.e.\ the
equivalence classes of diversified formulae) as objects of
$\ca^{st}$. We are going to establish a normalization procedure
for arrow terms of $\ca^{st}$ that eventually delivers our
coherence result. For this we use a sequence of definitions and
lemmata. We say that a form set $S$ is a \emph{subformset} of a
form set $T$ if there is a formula $A$ in $S$ and a formula $B$ in
$T$ ($S$ and $T$ are equivalence classes) such that $A$ is a
subformula of $B$. For example $p\kon(q\vee r)$ is a subformset of
$(r\vee q)\kon s\kon p$. We use freely for form sets the
terminology which is standard for formulae and say, for example,
that $r\vee q$ and $(r\vee q)\kon p$ are \emph{conjuncts} of the
form set $(r\vee q)\kon s\kon p$ whose \emph{main connective} is
$\kon$. We say that a conjunct $X$ of a form set is \emph{prime}
if $\kon$ is not the main connective in $X$. For example $r\vee q$
is a prime conjunct of $(r\vee q)\kon s\kon p$ but $(r\vee q)\kon
p$ is not. Also when $\kon$ is not the main connective of a form
set, we treat this form set as the prime conjunct of itself. We
use the same conventions for $\vee$ and, for example, $(r\vee
q)\kon s\kon p$ is the prime disjunct of itself. We denote by
$let(X)$ the set of letters in a form set $X$.

Every arrow term of $\ca^{st}$ is equal to a \emph{developed}
arrow term of the form
\[
s_n\cirk\ldots\cirk s_1\cirk\mj
\]
where every $s_i$ (if there is any) is a $\iota$-term. We tacitly
use developed form of arrow terms throughout the proofs of lemmata
given below. We take over the following lemma from \cite{DPInt}.

\prop{Lemma \thesection.1 (\cite{DPInt}, Section 14, Lemma~1)}{If
${u\!: X\str Y}$ is an arrow of $\ca^{st}$, and $P$ is a set of
letters such that for every subformset ${U\kon V}$ of $X$
\[ {let(U)}\subseteq
P\quad {\mbox {\it iff}} \quad{let(V)}\subseteq P,
\]
then this equivalence holds for every subformset ${U\kon V}$ of
$Y$.}

\noindent As a corollary (taking $P=let(X_1)$) we have the
following:

\prop{Lemma \thesection.2}{If $u\!:X_1\vee X_2\str Y$ is an arrow
of $\ca^{st}$, then for every subformset $U\kon V$ of $Y$ we have
that \[ let(U)\subseteq let(X_1)\quad\mbox{iff}\quad
let(V)\subseteq let(X_1).\] }

\vspace{-3ex}

\noindent (Since $X_1\vee X_2$ is the same form set as $X_2\vee
X_1$, it is not necessary to mention that the same holds when we
replace $X_1$ by $X_2$ in the conclusion of this lemma.)

\prop{Lemma \thesection.3}{Every arrow term $t\!:X'\kon X''\str Y$
of $\ca^{st}$ is equal to $t'\kon t''$ for some arrow terms
$t'\!:X'\str Y'$ and $t''\!:X''\str Y''$.}

\dkz Let $s_n\cirk\ldots\cirk s_1\cirk\mj$ be a developed arrow
term equal to $t$. We proceed by induction on $n$. If $n=0$, then
$t=\mj_{X'\kon X''}=\mj_{X'}\kon \mj_{X''}$. If $n>0$, by the
induction hypothesis we have $s_{n-1}\cirk\ldots\cirk
s_1\cirk\mj=u'\kon u''$ for some arrow terms $u'\!:X'\str Z'$ and
$u''\!:X''\str Z''$. So $s_n$ is a $\iota$-term whose source is
$Z'\kon Z''$. Since $\vee$ is the main connective of the source of
the head of $s_n$, we have that $s_n=v'\kon v''$ for some arrow
terms $v'\!:Z'\str Y'$ and $v''\!:Z''\str Y''$ where one of $v'$
and $v''$ is $\mj$ and the other is a $\iota$-term with the same
head as $s_n$. So, $t=(v'\kon v'')\cirk (u'\kon u'')=(v'\cirk
u')\kon (v''\cirk u'')$. \qed

Since the primitive equations of \ca\ and $\ca^{st}$ are such that
the number of occurrences of $\iota$ is the same on the both
sides, we have:

\prop{Lemma \thesection.4}{All the arrow terms representing the
same arrow of \ca\ or $\ca^{st}$ have the same number of
occurrences of $\iota$.}

We introduce now a procedure of deleting letters from form sets.
Roughly speaking, to delete a letter $p$ from a form set (which
includes some other letters) means to take a formula in this form
set, delete the letter $p$ together with its connective and
associated brackets from this formula, and then  form its
equivalence class. It is not difficult to see that this does not
depend on the choice of the formula in a form set. In terms of our
notation for form sets we define $X^{-p}$, for a form set $X$
different from $p$, in the following way:

\vspace{1ex}

if $p$ is not in $X$, then $X^{-p}$ is $X$;

\vspace{1ex}

if $X$ is of the form $Y\vee p$  or $Y\kon p$, then $X^{-p}$ is
$Y$;

\vspace{1ex}

if $X$ is of the form $Y\vee Z$ for $Y$ and $Z$ different from
$p$, then $X^{-p}$ is $Y^{-p}\vee Z^{-p}$, and the same holds when
we replace $\vee$ by $\kon$.

\vspace{2ex}

If $let(X)\setminus\{p,q\}\neq\emptyset$ then it is easy to see
that
\[
(X^{-p})^{-q}=(X^{-q})^{-p},
\]
and we can define, for a finite set $P=\{p_1,\ldots,p_n\}$ of
letters such that $let(X)\setminus P\neq\emptyset$,
\[
X^{-P}=_{df}(\mbox{\scriptsize
$\cdots$}(X^{-p_1})^{-p_2}\cdots)^{-p_n}.
\]
This can be extended to a procedure of letter deletion from the
arrow terms of $\ca^{st}$.

Let $u\!:X\str Y$ be an arrow term of $\ca^{st}$, and let $P$ be a
finite set of letters such that $let(X)\setminus P\neq\emptyset$
(hence $let(Y)\setminus P\neq\emptyset$, since $let(Y)=let(X)$)
and such that, as in Lemma \thesection.1, for every subformset
$U\kon V$ of $X$ we have $let(U)\subseteq P$ iff $let(V)\subseteq
P$. We define inductively the arrow term $u^{-P}\!:X^{-P}\str
Y^{-P}$ in the following way:

\vspace{1ex}

if $u$ is $\mj_X$, then $u^{-P}$ is $\mj_{X^{-P}}$;

\vspace{1ex}

if $u$ is $\iota_{[S,T,U,V]}$ then
\[
u^{-P}=_{df}\left\{
\begin{array}{l}
\mj_{X^{-P}},\quad \mbox{when}\hspace{1em} let(S\kon T)\subseteq
P\hspace{1em}\mbox{or}\hspace{1em}let(U\kon V)\subseteq P
\\[1ex]
\iota_{[S^{-P},T^{-P},U^{-P},V^{-P}]},\quad\mbox{otherwise;}
\end{array}
\right .
\]

\vspace{1ex}

if $u$ is $s\vee t$ for $s\!:S_1\str S_2$ and $t\!:T_1\str T_2$,
then
\[
u^{-P}=_{df}\left\{
\begin{array}{l}
s^{-P},\quad let(T_1)\subseteq P
\\[1ex]
t^{-P},\quad let(S_1)\subseteq P
\\[1ex]
s^{-P}\vee t^{-P},\quad\mbox{otherwise;}
\end{array}
\right .
\]
and we have the same clause when we replace $\vee$ by $\kon$;

\vspace{1ex}

if $u$ is $u_2\cirk u_1$, then by Lemma \thesection.1, both
$u_1^{-P}$ and $u_2^{-P}$ are defined and $u^{-P}$ is
$u_2^{-P}\cirk u_{1}^{-P}$.

\vspace{2ex}

Let $X_1$ and $X_2$ be form sets. We say that $\iota_{[S,T,U,V]}$
is $(X_1,X_2)$-\emph{splitting} when one of $let(S\kon T)$,
$let(U\kon V)$ is a subset of $let(X_1)$ while the other is a
subset of $let(X_2)$. We say that an arrow term of $\ca^{st}$ is
${(X_1,X_2)}$-\emph{splitting} when every occurrence of $\iota$ in
it is ${(X_1,X_2)}$-splitting, and we say that it is
${(X_1,X_2)}$-\emph{nonsplitting} when every occurrence of $\iota$
in it \emph{is not} ${(X_1,X_2)}$-splitting. For example,
$(\iota_{[p,q,s,t]}\kon\mj_{r\vee u})\cirk\iota_{[p\kon q,r,s\kon
t,u]}$ is a $(p\kon q\kon r,s\kon t\kon u)$-splitting arrow term.

One can easily check that if ${f=g}$ and $f$ is
${(X_1,X_2)}$-splitting, then $g$ is ${(X_1,X_2)}$-splitting, too.
This is not the case when we replace ``splitting'' by
``nonsplitting''. (Take for example the diagram $(7s)$ of the
preceding section and let $X_1$ be $U\kon X$ and $X_2$ be $V\kon
Y$, then the left leg of this diagram is $(X_1,X_2)$-nonsplitting,
and the occurrence of $\iota$ in $\iota\vee\mj$, in the right leg
is $(X_1,X_2)$-splitting.) It is clear that every
${(X_1,X_2)}$-splitting arrow term is equal to a developed
${(X_1,X_2)}$-splitting arrow term, and analogously with
``splitting'' replaced by ``nonsplitting''. We take over the
following three lemmata from \cite{DPInt}.

\prop{Lemma \thesection.5 (\cite{DPInt}, Section 14, Lemma~5)}{If
$u\!:X_1\vee X_2\str Y$ is ${(X_1,X_2)}$-nonsplitting, then $u$ is
equal to ${u_1\vee u_2}$ for some arrow terms $u_1\!:X_1\str X'_1$
and $u_2\!:X_2\str X'_2$.}

\noindent Note that for $X_1$, $X_2$, $X'_1$ and $X'_2$ as in
Lemma \thesection.5, an arrow term is $(X_1,X_2)$-splitting if and
only if it is $(X'_1,X'_2)$-splitting, which we will use later on.

\prop{Lemma \thesection.6 (\cite{DPInt}, Section 14, Lemma~6)}{If
$u\!:X_1\vee X_2\str Y$ is $(X_1,X_2)$-splitting, then $Y^{-X_1}$
is $X_2$ and $Y^{-X_2}$ is $X_1$.}

\prop{Lemma \thesection.7 (\cite{DPInt}, Section 14, Lemma~7)}{If
$u\!:X_1\vee X_2\str Y_1\kon Y_2$ is ${(X_1,X_2)}$-splitting, then
the main connective in $X_1$ and $X_2$ is $\kon$.}

Let $X_1=S\kon T$ and $X_2=U\kon V$ and let
$u\cirk\iota_{[S,T,U,V]}\!:X_1\vee X_2\str Y$ be
$(X_1,X_2)$-splitting. By Lemma \thesection.3, the main connective
in $Y$ is $\kon$ and by Lemma \thesection.2, the deletions
$^{-X_1}$ and $^{-X_2}$ are defined for every conjunct of Y. By
Lemma \thesection.6, we have $Y^{-X_2}=X_1=S\kon T$ and hence $Y$
is of the form $Y_S\kon Y_T$ for $Y_S$ and $Y_T$ such that
$Y_S^{-X_2}=S$ and $Y_T^{-X_2}=T$. Analogously, since
$Y^{-X_1}=X_2=U\kon V$, we have $Y=Y_U\kon Y_V$ for $Y_U^{-X_1}=U$
and $Y_V^{-X_1}=V$. We can then prove the following.

\prop{Lemma \thesection.8}{For $u\cirk\iota_{[S,T,U,V]}$ as above,
we have $Y_S=Y_U$ and $Y_T=Y_V$.}

\dkz Suppose $Y_S=Y_U\kon Z$, and hence, $Y_V=Y_T\kon Z$. We have
\[
\iota_{[S,T,U,V]}\!\!:\!(S\kon Y_T^{-X_2})\vee(U\kon
Y_T^{-X_1}\!\kon Z^{-X_1})\str(S\vee
U)\kon(Y_T^{-X_2}\vee(Y_T^{-X_1}\!\kon Z^{-X_1})).
\]
By Lemma \thesection.3, $u$ is of the form $s\kon t$ for
$t\!:Y_T^{-X_2}\vee(Y_T^{-X_1}\kon Z^{-X_1})\str W$, where $W$ is
a conjunct of $Y$. Since the source and the target of $t$ share
the same letters we have that $let(W)=let(Y_T)\cup let(Z^{-X_1})$.
Hence $W$ is of the form $Y_T\kon W'$ for $W'$ such that
$let(W')=let(Z^{-X_1})\subseteq X_2$. Since $W'$ is a conjunct of
$Y$, by Lemma \thesection.2, we have $let(Y)\subseteq X_2$ which
means that $let(X_1)=\emptyset$, i.e.\ a contradiction. We proceed
in the other cases quite similar. \qed

\vspace{2ex}

In the sequel, for $u\cirk\iota_{[S,T,U,V]}$ as above, we denote
by $Y_{SU}$ both $Y_S$ and $Y_U$ (which are equal by the preceding
lemma), and by the same reason we denote by $Y_{TV}$ both $Y_T$
and $Y_V$.

\prop{Lemma \thesection.9}{If $u\!:X_1\vee X_2\str Y_1\kon Y_2$ is
$(X_1,X_2)$-splitting, then $u$ factors as:}

\begin{center}
\begin{picture}(150,40)
\put(200,40){\makebox(0,0){\small $Y_1\kon Y_2$}}
\put(95,6){\makebox(0,0){\small $(Y_1^{-X_2}\vee
Y_1^{-X_1})\kon(Y_2^{-X_2}\vee Y_2^{-X_1})$}}
\put(0,40){\makebox(0,0){\small $(Y_1^{-X_2}\kon
Y_2^{-X_2})\vee(Y_1^{-X_1}\kon Y_2^{-X_1})$}}

\put(75,40){\vector(1,0){100}} \put(0,32){\vector(1,-1){20}}
\put(167,12){\vector(1,1){20}}

\put(185,22){\makebox(0,0)[l]{\small $u_1\kon u_2$}}
\put(5,22){\makebox(0,0)[r]{\small $\iota$}}
\put(125,44){\makebox(0,0)[b]{\small $u$}}

\end{picture}
\end{center}
{\it where $u_1\!:Y_1^{-X_2}\vee Y_1^{-X_1}\str Y_1$ is
$(Y_1^{-X_2},Y_1^{-X_1})$-splitting and $u_2\!:Y_2^{-X_2}\vee
Y_2^{-X_1}\str Y_2$ is $(Y_2^{-X_2},Y_2^{-X_1})$-splitting.}

\vspace{2ex}

\dkz We proceed by induction on number $n\geq 1$ of occurrences of
$\iota$ in $u$. First we prepare a ground for this induction. By
relying on the remark after the definition of
$(X_1,X_2)$-splitting arrow term and on Lemma \thesection.7, $u$
is equal to an arrow term of the form $v\cirk\iota_{[S,T,U,V]}$
for $X_1=S\kon T$, $X_2=U\kon V$, and for $v$, which by Lemma
\thesection.4 has $n-1$ occurrences of $\iota$, being
$(X_1,X_2)$-splitting.

If we denote $Y_1\kon Y_2$ by $Y$, then by Lemma \thesection.8 we
have $Y=Y_{SU}\kon Y_{TU}$ such that $Y_{SU}^{-X_2}=S$,
$Y_{SU}^{-X_1}=U$, $Y_{TV}^{-X_2}=T$ and $Y_{TV}^{-X_1}=V$. There
are several possibilities how the ``partition'' of $Y$ given by
the conjuncts $Y_1$ and $Y_2$ may be related to the one given by
the conjuncts $Y_{SU}$ and $Y_{TV}$, among which the following
three cases make all the essentially different situations.

\vspace{1ex}

(0) $\{Y_1,Y_2\}=\{Y_{SU},Y_{TV}\}$, when we are done;

\vspace{1ex}

(1) $Y_1\kon Z=Y_{SU}$ and $Y_2=Y_{TV}\kon Z$;

\vspace{1ex}

(2) $Y_1=Z_1\kon W_1$, $Y_2=Z_2\kon W_2$, $Y_{SU}=Z_1\kon Z_2$ and
$Y_{TV}=W_1\kon W_2$.

\vspace{1ex} \noindent We can start now with the induction.

If $n=1$, then $v=\mj_{(S\vee U)\kon(T\vee V)}$ and $\{S\vee
U,T\vee V\}=\{Y_1,Y_2\}$. Hence we are in case (0) and we are
done.

For the induction step suppose $n>1$. If we are in case (0), then
we are done.

If we are in case (1), then by Lemma \thesection.3 we have
$v=s\kon t$ for $(X_1,X_2)$-splitting arrow terms $s\!:S\vee U\str
Y_1\kon Z$ and $t\!:T\vee V\str Y_{TV}$ with less than $n$
occurrences of $\iota$ in them. Since $let(S)\subseteq let(X_1)$
and $let(U)\subseteq let(X_2)$, we have that $s$ is
$(S,U)$-splitting and we may apply the induction hypothesis in
order to obtain that $s$ is equal to an $(S,U)$-splitting (and
hence $(X_1,X_2)$-splitting) arrow term of the form

\begin{center}
\begin{picture}(215,20)
\put(260,10){\makebox(0,0){\small $Y_1\kon Z$.}}
\put(140,10){\makebox(0,0){\small $(Y_1^{-U}\vee
Y_1^{-S})\kon(Z^{-U}\vee Z^{-S})$}}
\put(0,10){\makebox(0,0){\small $(Y_1^{-U}\kon
Z^{-U})\vee(Y_1^{-S}\kon Z^{-S})$}}

\put(63,10){\vector(1,0){14}} \put(205,10){\vector(1,0){36}}

\put(223,14){\makebox(0,0)[b]{\small $s_1\kon s_2$}}
\put(70,14){\makebox(0,0)[b]{\small $\iota$}}

\end{picture}
\end{center}

Consider the following commutative diagram whose upper part is an
instance of $(7s)$:

\begin{center}
\begin{picture}(200,130)
\put(100,10){\makebox(0,0){\small $Y_1\kon Z\kon Y_{TV}=Y_1\kon
Y_2$}} \put(100,40){\makebox(0,0){\small $(Y_1^{-U}\vee
Y_1^{-S})\kon(Z^{-U}\vee Z^{-S})\kon(T\vee V)$}}

\put(20,90){\makebox(0,0){\small $((Y_1^{-U}\kon
Z^{-U})\vee(Y_1^{-S}\kon Z^{-S}))\kon(T\vee V)$}}
\put(180,70){\makebox(0,0){\small $(Y_1^{-U}\vee
Y_1^{-S})\kon((Z^{-U}\kon T)\vee(Z^{-S}\kon V))$}}
\put(100,120){\makebox(0,0){\small $(Y_1^{-U}\kon Z^{-U}\kon
T)\vee(Y_1^{-S}\kon Z^{-S}\kon V)=X_1\vee X_2$}}

\put(60,112){\vector(-4,-1){50}} \put(130,112){\vector(2,-1){60}}
\put(10,82){\vector(2,-1){60}} \put(190,62){\vector(-4,-1){50}}
\put(100,32){\vector(0,-1){14}}

\put(95,27){\makebox(0,0)[r]{\footnotesize $s_1\kon s_2\kon t$}}
\put(28,67){\makebox(0,0)[r]{\footnotesize $\iota\kon\mj$}}
\put(25,110){\makebox(0,0)[r]{\footnotesize $\iota$}}
\put(180,55){\makebox(0,0)[l]{\footnotesize $\mj\kon\iota$}}
\put(170,100){\makebox(0,0)[l]{\footnotesize $\iota$}}

\end{picture}
\end{center}
Since $Y_1^{-U}\kon Z^{-U}=S$ and $Y_1^{-S}\kon Z^{-S}=U$, the
left leg of this diagram is equal to $u$. Also, we have
$Y_1^{-U}=Y_1^{-X_2}$, $Z^{-U}=Z^{-X_2}$ (hence $Z^{-U}\kon
T=Y_2^{-X_2}$), $Y_1^{-S}=Y_1^{-X_1}$, and $Z^{-S}=Z^{-X_1}$
(hence $Z^{-S}\kon V=Y_2^{-X_1}$). So, the right leg of this
diagram is in the desired form.

If we are in case (2), then we use the induction hypothesis twice
and appeal to an instance of the following commutative diagram of
$\ca^{st}$ obtained by pasting instances of $(7s)$:

\begin{center}
\begin{picture}(200,90)
\put(100,5){\makebox(0,0){\small $(S\vee W)\kon(T\vee X)\kon(U\vee
Y)\kon(V\vee Z)$}}

\put(20,55){\makebox(0,0){\small $((S\kon T)\vee(W\kon
X))\kon((U\kon V)\vee(Y\kon Z))$}}
\put(180,35){\makebox(0,0){\small $((S\kon U)\vee(W\kon
Y))\kon((T\kon V)\vee(X\kon Z))$}}

\put(100,85){\makebox(0,0){$(S\kon T\kon U\kon V)\vee(W\kon X\kon
Y\kon Z)$}}

\put(25,40){\vector(2,-1){50}} \put(175,27){\vector(-4,-1){50}}
\put(75,77){\vector(-4,-1){50}} \put(125,77){\vector(2,-1){50}}

\put(38,28){\makebox(0,0)[r]{\footnotesize $\iota\kon\iota$}}
\put(35,74){\makebox(0,0)[r]{\footnotesize $\iota$}}
\put(165,18){\makebox(0,0)[l]{\footnotesize $\iota\kon\iota$}}
\put(150,71){\makebox(0,0)[l]{\footnotesize $\iota$}}

\put(260,0){\makebox(0,0)[l]{$\dashv$}}

\end{picture}
\end{center}

We have also the following three lemmata.

\prop{Lemma \thesection.10}{If $u\!:X_1\vee X_2\str Y'\vee Y''$ is
$(X_1,X_2)$-splitting and $Y'$ is a prime disjunct of $Y$ (i.e.\
$\vee$ is not the main connective in $Y'$), then $u=u'\vee u''$
for $u'\!:X'\str Y'$, where either\\[1ex] \indent for $i=1$ or $i=2$,
$X'$ is a prime disjunct of $X_i$ and $u'=\mj_{X'}$, or\\[1ex] \indent
$X'=X_1'\vee X_2'$ for $X_1'$ and $X_2'$ being prime disjuncts of
$X_1$ and $X_2$ respectively and $u'$ is $(X_1',X_2')$-splitting.}

\dkz By the dual of Lemma \thesection.3, $u$ is equal to $u'\vee
u''$, for $u'$ having $Y'$ as the target. If the source $X'$ of
$u'$ is a disjunct of $X_1$, by the assumption that $u$ is
$(X_1,X_2)$-splitting there are no occurrences of $\iota$ in $u'$
and hence it must be $\mj_{Y'}$, and $X'$, which is equal to $Y'$,
must be a prime disjunct of $X_1$.

If the source $X'$ of $u'$ is of the form $X_1'\vee X_2'$ for
$X_1'$ a disjunct of $X_1$ and $X_2'$ a disjunct of $X_2$, then
$Y'$ cannot be a letter and hence its main connective is $\kon$.
Also, $u'$ is $(X_1',X_2')$-splitting and by Lemma \thesection.7,
$X_1'$ and $X_2'$ are prime disjuncts of $X_1$ and $X_2$,
respectively. \qed

\prop{Lemma \thesection.11}{Let $v\cirk u\!:X_1\vee X_2\str Z$ be
such that $v$ is a $\iota$-term that is not $(X_1,X_2)$-splitting
and $u$ is an $(X_1,X_2)$-splitting arrow term. Then there exist
an arrow term $w$ and a $\iota$-term $v'$, which is not
$(X_1,X_2)$-splitting, such that $v\cirk u=w\cirk v'$. }

\dkz Let $X_1\vee X_2\stackrel{u\;}{\str} Y\stackrel{v\,}{\str}
Z$, and let $\iota_{[S,T,U,V]}$ be the head of $v$. We proceed by
induction on ``depth'' of $(S\kon T)\vee(U\kon V)$ in $Y$.

For the base of this induction we have the case when $S\kon T$ and
$U\kon V$ are prime disjuncts of $Y$. If $Y=(S\kon T)\vee(U\kon
V)\vee Y'''$, then by Lemma \thesection.10 we have
\[
v\cirk u=(\iota_{[S,T,U,V]}\cirk(u'\vee u''))\vee u'''
\]
for $u'\!:X'\str S\kon T$ and $u''\!:X''\str U\kon V$ satisfying
the conditions given by that lemma. (The arrow term
$u'''\!:X'''\str Y'''$ is out of our interest and it does not
exist when $Y=(S\kon T)\vee(U\kon V)$.)

We have several different situations depending on whether $X'$ or
$X''$ are prime disjuncts of $X_1$ or of $X_2$ or they are of the
form $X_1'\vee X_2'$ or $X_1''\vee X_2''$ for $X_1'$, $X_1''$
prime disjuncts of $X_1$ and $X_2'$, $X_2''$ prime disjuncts of
$X_2$. The following three cases represent essentially different
situations:

\vspace{1ex}

(0) For $i=1$ or $i=2$, $X'$ and $X''$ are prime disjuncts of
$X_i$. (By the assumption that $v$ is not $(X_1,X_2)$-splitting,
$X'$ and $X''$ cannot be prime disjuncts one of $X_1$ and the
other of $X_2$.) By Lemma \thesection.10, $u'$ and $u''$ are
identities and we are done.

\vspace{1ex}

(1) $X'=X_1'\vee X_2'$ for $X_1'$ and $X_2'$ being prime disjuncts
of $X_1$ and $X_2$ respectively, and $X''$ is a prime disjunct of
$X_1$. By Lemma \thesection.10, $u''=\mj_{X''}=\mj_{U\kon V}$ and
we may apply Lemma \thesection.9 to $u'\!:X_1'\vee X_2'\str S\kon
T$ which is $(X_1',X_2')$-splitting. So,
$\iota_{[S,T,U,V]}\cirk(u'\vee u'')$ is equal to the left leg of
the following commutative diagram whose upper part is an instance
of $(1s)$ and whose lower part is a naturality diagram for
$\iota$.

\begin{center}
\begin{picture}(200,140)
\put(120,10){\makebox(0,0){\small $(S\vee U)\kon(T\vee V)$}}
\put(20,40){\makebox(0,0){\small $(S\kon T)\vee(U\kon V)$}}

\put(120,60){\makebox(0,0){\small $(S^{-X_2'}\vee S^{-X_1'}\vee
U)\kon(T^{-X_2'}\vee T^{-X_1'}\vee V)$}}
\put(20,90){\makebox(0,0){\small $((S^{-X_2'}\vee
S^{-X_1'})\kon(T^{-X_2'}\vee T^{-X_1'}))\vee(U\kon V)$}}

\put(180,110){\makebox(0,0){\small $((S^{-X_2'}\vee
U)\kon(T^{-X_2'}\vee V))\vee(S^{-X_1'}\kon T^{-X_1'})$}}
\put(120,140){\makebox(0,0){\small $(S^{-X_2'}\kon
T^{-X_2'})\vee(S^{-X_1'}\kon T^{-X_1'})\vee(U\kon V)$}}

\put(35,34){\vector(4,-1){65}} \put(35,84){\vector(4,-1){65}}
\put(19,84){\vector(0,-1){38}} \put(119.5,54){\vector(0,-1){38}}
\put(135,134){\vector(4,-1){65}}
\put(98.5,130){\vector(-3,-1){79}}
\put(220,97){\vector(-3,-1){79}}

\put(190,83){\makebox(0,0)[l]{\footnotesize $\iota$}}
\put(75,80){\makebox(0,0)[l]{\footnotesize $\iota$}}
\put(65,20){\makebox(0,0)[l]{\footnotesize $\iota$}}
\put(15,65){\makebox(0,0)[r]{\footnotesize $(u_1'\kon
u_2')\vee\mj$}} \put(125,35){\makebox(0,0)[l]{\footnotesize
$(u_1'\vee\mj)\kon (u_2'\vee\mj)$}}
\put(45,120){\makebox(0,0)[r]{\footnotesize $\iota\vee\mj$}}
\put(185,128){\makebox(0,0)[l]{\footnotesize $\mj\vee\iota$}}

\end{picture}
\end{center}
The right leg of this diagram is of the desired form since it
starts with $\mj\vee\iota$ which is not $(X_1,X_2)$-splitting.

\vspace{1ex}

(2) $X'=X_1'\vee X_2'$ and $X''=X_1''\vee X_2''$ for $X_1'$ and
$X_1''$ being prime disjuncts of $X_1$, and $X_2'$ and $X_2''$
being prime disjuncts of $X_2$. Then we apply Lemma \thesection.9
to $u'\!:X_1'\vee X_2'\str S\kon T$ and to $u''\!:X_1''\vee
X_2''\str U\kon V$, and proceed as in case (1) relying on the
following commutative diagram of $\ca^{st}$ obtained by pasting
instances of $(1s)$:

\begin{center}
\begin{picture}(200,100)
\put(100,90){\makebox(0,0){\small $(S\kon W)\vee(T\kon
X)\vee(U\kon Y)\vee(V\kon Z)$}}

\put(20,60){\makebox(0,0){\small $((S\vee T)\kon(W\vee
X))\vee((U\vee V)\kon(Y\vee Z))$}}
\put(180,40){\makebox(0,0){\small $((S\vee U)\kon(W\vee
Y))\vee((T\vee V)\kon(X\vee Z))$}}

\put(100,10){\makebox(0,0){\small $(S\vee T\vee U\vee V)\kon(W\vee
X\vee Y\vee Z)$}}

\put(25,45){\vector(2,-1){50}} \put(175,32){\vector(-4,-1){50}}
\put(75,82){\vector(-4,-1){50}} \put(125,82){\vector(2,-1){50}}

\put(35,79){\makebox(0,0)[r]{\footnotesize $\iota\vee\iota$}}
\put(38,33){\makebox(0,0)[r]{\footnotesize $\iota$}}
\put(150,76){\makebox(0,0)[l]{\footnotesize $\iota\vee\iota$}}
\put(165,23){\makebox(0,0)[l]{\footnotesize $\iota$}}

\end{picture}
\end{center}

For the induction step, we proceed as follows. If $Y$ is of the
form $Y_1\kon Y_2$ where $Y_1$ is a prime conjunct of $Y$ whose
subformset is $(S\kon T)\vee(U\kon V)$, then, by Lemma
\thesection.9, $v\cirk u$ factors as:

\begin{center}
\begin{picture}(120,50)
\put(160,40){\makebox(0,0){\small $Z_1\kon Y_2$}}
\put(160,10){\makebox(0,0){\small $Y_1\kon Y_2$}}
\put(0,10){\makebox(0,0){\small $(Y_1^{-X_2}\vee
Y_1^{-X_1})\kon(Y_2^{-X_2}\vee Y_2^{-X_1})$}}
\put(0,40){\makebox(0,0){\small $(Y_1^{-X_2}\kon
Y_2^{-X_2})\vee(Y_1^{-X_1}\kon Y_2^{-X_1})$}}

\put(78,40){\vector(1,0){60}} \put(0,35){\vector(0,-1){20}}
\put(78,10){\vector(1,0){60}} \put(160,15){\vector(0,1){20}}

\put(165,25){\makebox(0,0)[l]{\small $v_1\kon \mj_{Y_2}$}}
\put(105,0){\makebox(0,0)[b]{\small $u_1\kon u_2$}}
\put(105,45){\makebox(0,0)[b]{\small $v\cirk u$}}
\put(-5,25){\makebox(0,0)[r]{\small $\iota$}}

\end{picture}
\end{center}
where $u_1$ is $(Y_1^{-X_2},Y_1^{-X_1})$-splitting and $v_1$ is a
$\iota$-term that is not $(X_1,X_2)$-splitting, and hence it is
not $(Y_1^{-X_2},Y_1^{-X_1})$-splitting. By the induction
hypothesis $v_1\cirk u_1$ is equal to an arrow term of the form
$w_1\cirk v_1'$ for $v_1'$ a $\iota$-term that is not
$(Y_1^{-X_2},Y_1^{-X_1})$-splitting. By Lemma \thesection.5, we
may assume that $v_1'$ is of the form $v_1''\vee\mj_{Y_1^{-X_1}}$
or $\mj_{Y_1^{-X_2}}\vee v_1''$. In both cases we just apply the
naturality of $\iota$ and we are done.

If $Y$ is of the form $Y'\vee Y''$ where $Y'$ is the prime
disjunct containing $(S\kon T)\vee(U\kon V)$ as a subformset, then
by Lemma \thesection.10, $v\cirk u$ is of the form

\begin{center}
\begin{picture}(215,20)
\put(200,10){\makebox(0,0){\small $Z'\vee Y''$,}}
\put(100,10){\makebox(0,0){\small $Y'\vee Y''$}}
\put(0,10){\makebox(0,0){\small $X'\vee X''$}}

\put(30,10){\vector(1,0){40}} \put(130,10){\vector(1,0){40}}

\put(50,14){\makebox(0,0)[b]{\small $u'\vee u''$}}
\put(150,14){\makebox(0,0)[b]{\small $v'\vee\mj_{Y''}$}}

\end{picture}
\end{center}
and if for $i=1$ or $i=2$, $X'$ is a prime disjunct of $X_i$, then
$u'$ is $\mj_{X'}$ and we are done. If $X'=X_1'\vee X_2'$, then we
just apply the induction hypothesis to $v'\cirk u'$.
\mbox{\hspace{1em}}\qed

\prop{Lemma \thesection.12}{For every arrow term $t\!:X_1\vee
X_2\str Y$ there are arrow terms $v_1\!:X_1\str Y^{-X_2}$,
$v_2\!:X_2\str Y^{-X_1}$ and an $(X_1,X_2)$-splitting arrow term
$u\!:Y^{-X_2}\vee Y^{-X_1}\str Y$ such that $t=u\cirk(v_1\vee
v_2)$.}

\dkz We proceed by induction on the number $n$ of occurrences of
$\iota$ in $t$. If $n=0$, then since identities are at the same
time $(X_1,X_2)$-splitting and $(X_1,X_2)$-nonsplitting we are
done.

For the induction step, take a developed arrow term equal to $t$.
If every $\iota$-term in it is $(X_1,X_2)$-splitting, then by
Lemma \thesection.6 we are done. Otherwise, by Lemma
\thesection.11 (applied to this developed arrow term from its
right-hand side end up to the rightmost $\iota$-term in it that is
not $(X_1,X_2)$-splitting) we have $t=t'\cirk v'$ where $v'$ is
not $(X_1,X_2)$-splitting $\iota$-term. By Lemma \thesection.5,
$v'=v_1'\vee v_2'$ for $v_1'\!:X_1\str X_1'$ and $v_2'\!:X_2\str
X_2'$. By Lemma \thesection.4, $t'$ has $n-1$ occurrences of
$\iota$ and since $let(X_1)=let(X_1')$ and $let(X_2)=let(X_2')$,
we may apply the induction hypothesis to it. \qed

\vspace{2ex}

We conclude this section with the following proof.

\vspace{2ex}

\noindent {\sc Proof of Proposition~4.2.}\quad Let $t\!:X\str Y$
be an arrow of $\ca^{st}$. To prove that $t$ is unique, we proceed
by induction on the complexity of $X$ and $Y$. If $X$ is a letter
$p$, then $Y$ must be $p$ too, and $t\!:p\str p$ must be $\mj_p$.

If $X=X'\kon X''$, then by Lemma \thesection.3 and the induction
hypothesis, $t=t'\kon t''$ for unique arrows $t'\!:X'\str Y'$ and
$t''\!:X''\str Y''$. We reason analogously when $Y=Y'\vee Y''$.

Suppose $X=X_1\vee X_2$ and $Y=Y_1\kon Y_2$. Then by Lemmata
\thesection.12 and \thesection.9, and the induction hypothesis,
$t$ is equal to the following composition

\begin{center}
\begin{picture}(310,25)
\put(310,10){\makebox(0,0){\small $Y_1\kon Y_2$}}
\put(195,10){\makebox(0,0){\small $(Y_1^{-X_2}\vee
Y_1^{-X_1})\kon(Y_2^{-X_2}\vee Y_2^{-X_1})$}}
\put(0,10){\makebox(0,0){\small $X_1\vee X_2$}}
\put(80,10){\makebox(0,0){\small $Y^{-X_2}\vee Y^{-X_1}$}}

\put(20,10){\vector(1,0){27}} \put(267,10){\vector(1,0){27}}
\put(110,10){\vector(1,0){12}}

\put(280,16){\makebox(0,0)[b]{\small $u_1\kon u_2$}}
\put(118,16){\makebox(0,0)[b]{\small $\iota$}}
\put(33,16){\makebox(0,0)[b]{\small $v_1\vee v_2$}}

\end{picture}
\end{center}
for unique arrows $u_1$, $u_2$, $v_1$ and $v_2$. (Note that all
the sources and the targets above are completely determined by
$X_1$, $X_2$, $Y_1$ and $Y_2$.) So, $t$ is the unique arrow with
the source $X$ and the target $Y$. \qed

\section{A note on reduced bar construction}

This section is optional. Its aim is to give an analysis of a
reduced bar construction based on a monoid in a category whose
monoidal structure is given by finite products. Such a reduced bar
construction was used by Thomason in \cite{T79}. We believe this
analysis is not new, but we couldn't find (or just couldn't
recognize) a reference which covers it completely, especially in
its graphical approach we intend to use.

Let $\Delta$ be algebraist's simplicial category defined as in
\cite{ML71}, VII.5, for whose arrows we take over the notation
used in that book. Let $\Delta_+$ (denoted by $\Delta^+$ in
\cite{ML71}) be the topologist's simplicial category which is the
full subcategory of $\Delta$ with objects all nonempty ordinals
$\{1,2,3,\ldots\}$. In order to use geometric dimension, the
objects of $\Delta_+$ are rewritten as $\{0,1,2,\ldots\}$. So,
$\Delta_+$ has all finite ordinals as objects, and in this
category the source of $\delta_i^n$, for $n\geq 1$ and $0\leq
i\leq n$, is $n-1$ and the target is $n$, while the source of
$\sigma_i^n$, for $n\geq 1$ and $0\leq i\leq n-1$ is $n$ and the
target is $n-1$. When we speak of $\Delta_+^{op}$, then we denote
its arrows $(\delta_i^n)^{op}\!:n\str n-1$ by $d_i^n$ and
$(\sigma_i^n)^{op}\!:n-1\str n$ by $s_i^n$.

It is known that the functor ${\cal J}\!:\Delta^{op}\str\Delta$
defined on objects as ${\cal J}(n)=n\pl 1$, and on arrows by the
clauses

\vspace{1ex}

\begin{center}
\begin{picture}(320,40)

\put(0,20){\makebox(0,0)[l]{${\cal
J}(\delta_i^n)^{op}=\sigma_i^{n+1}$}}

\put(100,10){\circle*{2}} \put(130,10){\circle*{2}}
\put(145,10){\circle*{2}} \put(170,10){\circle*{2}}
\put(100,30){\circle*{2}} \put(130,30){\circle*{2}}
\put(145,30){\circle{2}} \put(160,30){\circle*{2}}
\put(185,30){\circle*{2}}

\put(100,0){\makebox(0,0)[b]{\scriptsize $0$}}
\put(130,0){\makebox(0,0)[b]{\scriptsize $i\mn 1$}}
\put(145,0){\makebox(0,0)[b]{\scriptsize $i$}}
\put(170,0){\makebox(0,0)[b]{\scriptsize $n\mn 1$}}

\put(100,35){\makebox(0,0)[b]{\scriptsize $0$}}
\put(130,35){\makebox(0,0)[b]{\scriptsize $i\mn 1$}}
\put(145,35){\makebox(0,0)[b]{\scriptsize $i$}}
\put(160,35){\makebox(0,0)[b]{\scriptsize $i\pl 1$}}
\put(185,35){\makebox(0,0)[b]{\scriptsize $n$}}

\put(100,29){\line(0,-1){18}} \put(130,29){\line(0,-1){18}}

\put(145.7,10.7){\line(3,4){14}} \put(170.7,10.7){\line(3,4){14}}

\put(116,20){\makebox(0,0){\ldots}}
\put(166,20){\makebox(0,0){\ldots}}

\put(200,20){\makebox(0,0)[l]{$\mapsto$}}

\put(220,10){\circle*{2}} \put(250,10){\circle*{2}}
\put(265,10){\circle*{2}} \put(280,10){\circle*{2}}
\put(305,10){\circle*{2}} \put(220,30){\circle*{2}}
\put(250,30){\circle*{2}} \put(265,30){\circle*{2}}
\put(280,30){\circle*{2}} \put(295,30){\circle*{2}}
\put(320,30){\circle*{2}}

\put(220,0){\makebox(0,0)[b]{\scriptsize $0$}}
\put(265,0){\makebox(0,0)[b]{\scriptsize $i$}}
\put(305,0){\makebox(0,0)[b]{\scriptsize $n$}}

\put(220,35){\makebox(0,0)[b]{\scriptsize $0$}}
\put(265,35){\makebox(0,0)[b]{\scriptsize $i$}}
\put(280,35){\makebox(0,0)[b]{\scriptsize $i\pl 1$}}
\put(320,35){\makebox(0,0)[b]{\scriptsize $n\pl 1$}}

\put(220,29){\line(0,-1){18}} \put(250,29){\line(0,-1){18}}
\put(265,29){\line(0,-1){18}}

\put(265.7,10.7){\line(3,4){14}} \put(280.7,10.7){\line(3,4){14}}
\put(305.7,10.7){\line(3,4){14}}

\put(236,20){\makebox(0,0){\ldots}}
\put(301,20){\makebox(0,0){\ldots}}

\end{picture}
\end{center}

\vspace{1ex}

\begin{center}
\begin{picture}(320,40)

\put(0,20){\makebox(0,0)[l]{${\cal
J}(\sigma_i^n)^{op}=\delta_{i+1}^{n+1}$}}

\put(100,30){\circle*{2}} \put(130,30){\circle*{2}}
\put(145,30){\circle*{2}} \put(170,30){\circle*{2}}
\put(100,10){\circle*{2}} \put(130,10){\circle*{2}}
\put(145,10){\circle*{2}} \put(160,10){\circle*{2}}
\put(185,10){\circle*{2}}

\put(100,35){\makebox(0,0)[b]{\scriptsize $0$}}
\put(130,35){\makebox(0,0)[b]{\scriptsize $i\mn 1$}}
\put(145,35){\makebox(0,0)[b]{\scriptsize $i$}}
\put(170,35){\makebox(0,0)[b]{\scriptsize $n\mn 1$}}

\put(100,0){\makebox(0,0)[b]{\scriptsize $0$}}
\put(130,0){\makebox(0,0)[b]{\scriptsize $i\mn 1$}}
\put(145,0){\makebox(0,0)[b]{\scriptsize $i$}}
\put(160,0){\makebox(0,0)[b]{\scriptsize $i\pl 1$}}
\put(185,0){\makebox(0,0)[b]{\scriptsize $n$}}

\put(100,29){\line(0,-1){18}} \put(130,29){\line(0,-1){18}}
\put(145,29){\line(0,-1){18}}

\put(145.7,29.3){\line(3,-4){14}}
\put(170.7,29.3){\line(3,-4){14}}
\put(116,20){\makebox(0,0){\ldots}}
\put(166,20){\makebox(0,0){\ldots}}

\put(200,20){\makebox(0,0)[l]{$\mapsto$}}

\put(220,30){\circle*{2}} \put(265,30){\circle*{2}}
\put(280,30){\circle*{2}} \put(305,30){\circle*{2}}
\put(220,10){\circle*{2}} \put(265,10){\circle*{2}}
\put(280,10){\circle{2}} \put(295,10){\circle*{2}}
\put(320,10){\circle*{2}}

\put(220,35){\makebox(0,0)[b]{\scriptsize $0$}}
\put(265,35){\makebox(0,0)[b]{\scriptsize $i$}}
\put(305,35){\makebox(0,0)[b]{\scriptsize $n$}}

\put(220,0){\makebox(0,0)[b]{\scriptsize $0$}}
\put(265,0){\makebox(0,0)[b]{\scriptsize $i$}}
\put(280,0){\makebox(0,0)[b]{\scriptsize $i\pl 1$}}
\put(320,0){\makebox(0,0)[b]{\scriptsize $n\pl 1$}}

\put(220,29){\line(0,-1){18}} \put(265,29){\line(0,-1){18}}

\put(280.7,29.3){\line(3,-4){14}}
\put(305.7,29.3){\line(3,-4){14}}

\put(244,20){\makebox(0,0){\ldots}}
\put(301,20){\makebox(0,0){\ldots}}

\end{picture}
\end{center}
is faithful and obviously injective on objects. The most explicit
formulation of this result is given in the fifth paragraph of
\cite{T07}, where the author refers to \cite{L06} as a paper in
which this idea figures prominently. (See also \cite{DP08a},
Section~6, and \cite{DP08b}, Section~6, for some more general
results). Intuitively, this functor is given by taking complements
of the standard graphical presentations for the arrows of
$\Delta^{op}$ as in the following picture, where the inner graph
with solid lines represents an arrow of $\Delta^{op}$, and the
outer graph with dotted lines represents its image under $\cal J$.

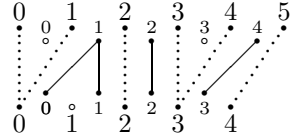
\begin{center}
\begin{picture}(100,50)

\put(10,15){\circle*{2}} \put(30,15){\circle*{2}}
\put(50,15){\circle*{2}} \put(70,15){\circle*{2}}

\put(10,35){\circle{2}} \put(30,35){\circle*{2}}
\put(50,35){\circle*{2}} \put(70,35){\circle{2}}
\put(90,35){\circle*{2}}

\put(0,10){\circle*{2}} \put(20,10){\circle{2}}
\put(40,10){\circle*{2}} \put(60,10){\circle*{2}}
\put(80,10){\circle*{2}}

\put(0,40){\circle*{2}} \put(20,40){\circle*{2}}
\put(40,40){\circle*{2}} \put(60,40){\circle*{2}}
\put(80,40){\circle*{2}} \put(100,40){\circle*{2}}

\put(10,12){\makebox(0,0)[t]{\scriptsize $0$}}
\put(30,12){\makebox(0,0)[t]{\scriptsize $1$}}
\put(50,12){\makebox(0,0)[t]{\scriptsize $2$}}
\put(70,12){\makebox(0,0)[t]{\scriptsize $3$}}
\put(10,12){\makebox(0,0)[t]{\scriptsize $0$}}

\put(10,38){\makebox(0,0)[b]{\scriptsize $0$}}
\put(30,38){\makebox(0,0)[b]{\scriptsize $1$}}
\put(50,38){\makebox(0,0)[b]{\scriptsize $2$}}
\put(70,38){\makebox(0,0)[b]{\scriptsize $3$}}
\put(90,38){\makebox(0,0)[b]{\scriptsize $4$}}

\put(0,7){\makebox(0,0)[t]{$0$}} \put(20,7){\makebox(0,0)[t]{$1$}}
\put(40,7){\makebox(0,0)[t]{$2$}}
\put(60,7){\makebox(0,0)[t]{$3$}}
\put(80,7){\makebox(0,0)[t]{$4$}}

\put(0,43){\makebox(0,0)[b]{$0$}}
\put(20,43){\makebox(0,0)[b]{$1$}}
\put(40,43){\makebox(0,0)[b]{$2$}}
\put(60,43){\makebox(0,0)[b]{$3$}}
\put(80,43){\makebox(0,0)[b]{$4$}}
\put(100,43){\makebox(0,0)[b]{$5$}}

\put(10,15){\line(1,1){20}} \put(30,15){\line(0,1){20}}
\put(50,15){\line(0,1){20}} \put(70,15){\line(1,1){20}}

\multiput(0,10)(0,3){10}{\circle*{.5}}
\multiput(0,10)(2,3){10}{\circle*{.5}}
\multiput(40,10)(0,3){10}{\circle*{.5}}
\multiput(60,10)(0,3){10}{\circle*{.5}}
\multiput(60,10)(2,3){10}{\circle*{.5}}
\multiput(80,10)(2,3){10}{\circle*{.5}}

\end{picture}
\end{center}

So, we may regard of $\Delta^{op}$ as a subcategory of $\Delta$.
From now on we restrict $\cal J$ to $\Delta_+^{op}$ taking into
account that, this time, it is defined on objects by the clause
${\cal J}(n)=n\pl 2$.

Let $\Delta_{Int}$ be the subcategory of $\Delta$ whose objects
are finite ordinals greater or equal to 2 and whose arrows are
interval maps, i.e.\ order-preserving functions, which preserve,
moreover, the first and the last element. The category
$\Delta_{Int}$ is the image  of $\Delta_+^{op}$ under the functor
$\cal J$. So, $\Delta_{Int}$ is isomorphic to $\Delta_+^{op}$ and
in the sequel we will represent the arrows of $\Delta_+^{op}$ by
the standard graphical presentations for the corresponding arrows
of $\Delta_{Int}$.

Let $\Delta_{par}$ be the category whose objects are again finite
ordinals and whose arrows are order preserving partial functions.
Beside the arrows $\delta_i^n$ and $\sigma_i^n$, to generate
$\Delta_{par}$ we need also the arrows $\rho_i^n\!:n\pl 1\str n$
for $n\geq 0$ and $0\leq i\leq n$, which are partial functions
graphically presented as

\begin{center}
\begin{picture}(85,40)

\put(0,10){\circle*{2}} \put(30,10){\circle*{2}}
\put(45,10){\circle*{2}} \put(70,10){\circle*{2}}
\put(0,30){\circle*{2}} \put(30,30){\circle*{2}}
\put(45,30){\circle{2}} \put(60,30){\circle*{2}}
\put(85,30){\circle*{2}}

\put(0,0){\makebox(0,0)[b]{\scriptsize $0$}}
\put(30,0){\makebox(0,0)[b]{\scriptsize $i\mn 1$}}
\put(45,0){\makebox(0,0)[b]{\scriptsize $i$}}
\put(70,0){\makebox(0,0)[b]{\scriptsize $n\mn 1$}}

\put(0,35){\makebox(0,0)[b]{\scriptsize $0$}}
\put(30,35){\makebox(0,0)[b]{\scriptsize $i\mn 1$}}
\put(45,35){\makebox(0,0)[b]{\scriptsize $i$}}
\put(60,35){\makebox(0,0)[b]{\scriptsize $i\pl 1$}}
\put(85,35){\makebox(0,0)[b]{\scriptsize $n$}}

\put(0,29){\line(0,-1){18}} \put(30,29){\line(0,-1){18}}

\put(45.7,10.7){\line(3,4){14}} \put(70.7,10.7){\line(3,4){14}}

\put(16,20){\makebox(0,0){\ldots}}
\put(66,20){\makebox(0,0){\ldots}}

\end{picture}
\end{center}
The standard list of equations that satisfy $\delta$'s and
$\sigma$'s should be extended by the following equations:
\[
\rho_j\rho_i=\rho_i\rho_{j+1}\mbox{\hspace{2.3em}} i\leq j
\]
\[
\rho_j\delta_i=\left\{
\begin{array}{ll}
\delta_{i-1}\rho_j & i>j
\\[1ex]
\mj & i=j
\\[1ex]
\delta_i\rho_{j-1} & i<j
\end{array}
\right . \quad\quad\quad \rho_j\sigma_i=\left\{
\begin{array}{ll}
\sigma_{i-1}\rho_j & i>j
\\[1ex]
\rho_i\rho_i & i=j
\\[1ex]
\sigma_i\rho_{j+1} & i<j
\end{array}
\right .
\]

A \emph{counital monad} $\langle T,\eta,\mu,\varepsilon\rangle$ in
a category $X$ consists of a functor $T\!:X\str X$ and three
natural transformations
\[
\eta\!:{\cal I}_X\strt T,\quad \mu\!:T^2\strt
T\quad\mbox{and}\quad\varepsilon\!:T\strt{\cal I}_X
\]
such that $\langle T,\eta,\mu\rangle$ is a monad in $X$, and
moreover,
\[
\varepsilon\cirk\eta=\mj_{{\cal
I}_X},\quad\varepsilon\cirk\mu=\varepsilon\cirk\varepsilon_T.
\]

%%%%%%%%%%%%%%%%%%%%%%%%%%%%%%%%%%%%%%%%%%%%%%%%%%%%%%%%%%%%%%%%%
In order to show that $\Delta_{par}$ bears a structure of a freely
generated counital monad, we construct an auxiliary syntactical
category $CM_0$. The construction of $CM_0$ is analogous to the
construction of \cm\ given at the beginning of Section~3, save
that the set $\cal P$ of generators is now replaced by a singleton
set $\{0\}$. The objects of $CM_0$ are finite ordinals, where $n$
stands for a sequence of $n$ occurrences of $T$; so $Tn$ is
${n+1}$. The arrows of $CM_0$ are defined syntactically as
equivalence classes of arrow terms generated from primitive arrow
terms $\mj_n\!:n\str n$, $\eta_n\!:n\str n+1$, $\mu_n\!:n+2\str
n+1$ and $\varepsilon_n\!:n+1\str n$, with the help of $\cirk$ and
$T$. These equivalence classes are taken with respect to the
smallest equivalence relation on arrow terms which makes out of
$\langle T,\eta,\mu,\varepsilon\rangle$ a counital monad in the
category $CM_0$.

The category $CM_0$ together with its counital monad structure is
freely generated in the following sense. It is the image of a
singleton set $\{0\}$ under the left adjoint of the forgetful
functor from the category of small counital monads (whose arrows
are functors preserving the counital monad structure on the nose,
i.e.\ exactly) to the category \emph{Set}; this forgetful functor
assigns to a counital monad the set of objects of its underlying
category. We can prove the following.

\prop{Proposition \thesection.1}{The categories $CM_0$ and
$\Delta_{par}$ are isomorphic.}

\dkz Consider the functor $G\!:CM_0\str\Delta_{par}$, which is
identity on objects and on arrows is defined so that
$G(T^m\eta_n)=\delta_m^{n+m}$, $G(T^m\mu_n)=\sigma_m^{m+n+1}$ and
$G(T^m\varepsilon_n)=\rho_m^{n+m}$. The equations of $CM_0$ enable
us to find for every arrow term $f$ of $CM_0$ an arrow term equal
to $f$ in the \emph{normal form} $f_4\cirk f_2\cirk f_1$ where
$f_1$ is free of $\mu$ and $\eta$, $f_2$ is free of $\varepsilon$
and $\eta$, and $f_4$ is free of $\varepsilon$ and $\mu$. Then we
may continue reasoning as in the proof of $S4_{\Box\Diamond}$
Coherence given in \cite{DP08b} (Section~4) to conclude that $G$
is faithful. (For this we replace the ordinals in the targets of
$f_1$, $f_2$ and $f_4$ by occurrences of $\diamond$'s.) So, we may
conclude that $G$ is an isomorphism. \mbox{\hspace{2em}}\qed

Let $\langle{\cal K},\otimes,I\rangle$ be a monoidal category. A
\emph{counital monoid} in $\cal K$ is a quadruple $\langle
C,\mu,\eta,\varepsilon\rangle$, where $\langle C,\mu,\eta\rangle$
is a monoid and $\varepsilon\!:C\str I$ is a monoid morphism, for
$I$ being equipped with its canonical monoid structure. The
category $\Delta_{par}$ is a strict monoidal category with $+$ as
tensor and $0$ as monoidal unit. The object $1$ of this category
together with $\sigma_0^1\!:2\str 1$ and $\delta_0^0\!:0\str 1$ is
a monoid equipped with the monoid morphism $\rho_0^0\!:1\str 0$.
So $\langle 1,\sigma_0^1,\delta_0^0,\rho_0^0\rangle$ is a counital
monoid. It is universal in the following sense (cf.\ Proposition~1
of \cite{ML71}, Section VII.5).

\prop{Proposition \thesection.2}{Given a counital monoid $\langle
C,\mu',\eta',\varepsilon'\rangle$ in a strict monoidal category
$\langle{\cal K},\otimes,I\rangle$, there is a unique strict
monoidal functor ${\cal F}\!:\Delta_{par}\str{\cal K}$ such that
${\cal F}(1)=C$, ${\cal F}(\sigma_0^1)=\mu'$, ${\cal
F}(\delta_0^0)=\eta'$ and ${\cal F}(\rho_0^0)=\varepsilon'$.}

\dkz The functor $C\otimes\underline{\hspace{.7em}}\!:{\cal
K}\str{\cal K}$ together with the natural transformations $\eta$,
$\mu$ and $\varepsilon$ derived, with the help of $\otimes$, from
the morphisms $\eta'$, $\mu'$ and $\varepsilon'$ respectively,
make a counital monad in $\cal K$. Since $CM_0$ together with its
counital monad structure is freely generated by $\{0\}$, by
relying on Proposition \thesection.1, we obtain a unique functor
${\cal F}\!:\Delta_{par}\str{\cal K}$, which maps the generator
$0$ to $I$ and preserves the counital monad structure of
$\Delta_{par}$ inherited from $CM_0$. This guarantees that ${\cal
F}$ is strict monoidal, and moreover, ${\cal F}(1)={\cal
F}(T0)=C\otimes I=C$, ${\cal F}(\sigma_0^1)={\cal
F}(G(\mu_0))=\mu'$, ${\cal F}(\delta_0^0)={\cal
F}(G(\eta_0))=\eta'$ and ${\cal F}(\rho_0^0)={\cal
F}(G(\varepsilon_0))=\varepsilon'$. {\mbox{\hspace{2em}}}\qed

%%%%%%%%%%%%%%%%%%%%%%%%%%%%%%%%%%%%%%%%%%%%%%%%%%%%%%%%%

We have that $\Delta_{Int}$ is a subcategory of $\Delta_{par}$ and
also we have a functor ${{\cal H}\!:\Delta_{Int}\str\Delta_{par}}$
defined on objects as ${\cal H}(n)=n\mn 2$, and on arrows, for
$f\!:n\str m$, as
\[
{\cal H}(f)=\rho_0^{m-2}\cirk\rho_{m-1}^{m-1}\cirk
f\cirk\delta_{n-1}^{n-1}\cirk\delta_0^{n-2}.
\]
(Intuitively, ${\cal H}(f)$ is obtained by omitting points $0$,
$n\mn 1$ from the source, and $0$, $m\mn 1$ from the target in the
graphical presentation of $f$ together with all edges including
them.) By using essentially the property that the arrows of
$\Delta_{Int}$ may be built free of $\delta_0^n$ and $\delta_n^n$,
it is not difficult to check that ${\cal H}$ so defined is indeed
a functor. (Note that ${\cal H}$ is not a functor from $\Delta$ to
$\Delta_{par}$.)

The composition ${\cal H}\cirk{\cal J}$ is a functor from
$\Delta_+^{op}$ to $\Delta_{par}$ which is identity on objects. In
this way $d_0^n$ and $d_n^n$ are mapped to the partial functions
graphically presented as

\begin{center}
\begin{picture}(215,40)

\put(-30,20){\makebox(0,0)[r]{$d_0^n$}}
\put(-15,20){\makebox(0,0){$\mapsto$}}

\put(0,10){\circle*{2}} \put(50,10){\circle*{2}}
\put(0,30){\circle{2}} \put(15,30){\circle*{2}}
\put(65,30){\circle*{2}}

\put(0,0){\makebox(0,0)[b]{\scriptsize $0$}}
\put(50,0){\makebox(0,0)[b]{\scriptsize $n\mn 2$}}

\put(0,35){\makebox(0,0)[b]{\scriptsize $0$}}
\put(15,35){\makebox(0,0)[b]{\scriptsize $1$}}
\put(65,35){\makebox(0,0)[b]{\scriptsize $n\mn 1$}}

\put(0.7,10.7){\line(3,4){14}} \put(50.7,10.7){\line(3,4){14}}

\put(35,20){\makebox(0,0){\ldots}}

\put(120,20){\makebox(0,0)[r]{$d_n^n$}}
\put(135,20){\makebox(0,0){$\mapsto$}}

\put(150,10){\circle*{2}} \put(200,10){\circle*{2}}
\put(150,30){\circle*{2}} \put(200,30){\circle*{2}}
\put(215,30){\circle{2}}

\put(150,0){\makebox(0,0)[b]{\scriptsize $0$}}
\put(200,0){\makebox(0,0)[b]{\scriptsize $n\mn 2$}}

\put(150,35){\makebox(0,0)[b]{\scriptsize $0$}}
\put(197,35){\makebox(0,0)[b]{\scriptsize $n\mn 2$}}
\put(218,35){\makebox(0,0)[b]{\scriptsize $n\mn 1$}}

\put(150,29){\line(0,-1){18}} \put(200,29){\line(0,-1){18}}

\put(175,20){\makebox(0,0){\ldots}}

\end{picture}
\end{center}
while $d_i^n$, for $0<i<n$, and $s_i^n$ are mapped to the
functions graphically presented as

\begin{center}
\begin{picture}(215,40)

\put(-30,20){\makebox(0,0)[r]{$d_i^n$}}
\put(-15,20){\makebox(0,0){$\mapsto$}}

\put(0,10){\circle*{2}} \put(30,10){\circle*{2}}
\put(60,10){\circle*{2}} \put(0,30){\circle*{2}}
\put(30,30){\circle*{2}} \put(45,30){\circle*{2}}
\put(75,30){\circle*{2}}

\put(0,0){\makebox(0,0)[b]{\scriptsize $0$}}
\put(30,0){\makebox(0,0)[b]{\scriptsize $i\mn 1$}}
\put(60,0){\makebox(0,0)[b]{\scriptsize $n\mn 2$}}

\put(0,35){\makebox(0,0)[b]{\scriptsize $0$}}
\put(30,35){\makebox(0,0)[b]{\scriptsize $i\mn 1$}}
\put(45,35){\makebox(0,0)[b]{\scriptsize $i$}}
\put(75,35){\makebox(0,0)[b]{\scriptsize $n\mn 1$}}

\put(0,29){\line(0,-1){18}} \put(30,29){\line(0,-1){18}}
\put(30.7,10.7){\line(3,4){14}} \put(60.7,10.7){\line(3,4){14}}

\put(15,20){\makebox(0,0){\ldots}}
\put(55,20){\makebox(0,0){\ldots}}

\put(120,20){\makebox(0,0)[r]{$s_i^n$}}
\put(135,20){\makebox(0,0){$\mapsto$}}

\put(150,10){\circle*{2}} \put(180,10){\circle*{2}}
\put(195,10){\circle{2}} \put(210,10){\circle*{2}}
\put(240,10){\circle*{2}} \put(150,30){\circle*{2}}
\put(180,30){\circle*{2}} \put(195,30){\circle*{2}}
\put(225,30){\circle*{2}}

\put(150,0){\makebox(0,0)[b]{\scriptsize $0$}}
\put(180,0){\makebox(0,0)[b]{\scriptsize $i\mn 1$}}
\put(195,0){\makebox(0,0)[b]{\scriptsize $i$}}
\put(210,0){\makebox(0,0)[b]{\scriptsize $i\pl 1$}}
\put(240,0){\makebox(0,0)[b]{\scriptsize $n\mn 1$}}

\put(150,35){\makebox(0,0)[b]{\scriptsize $0$}}
\put(180,35){\makebox(0,0)[b]{\scriptsize $i\mn 1$}}
\put(195,35){\makebox(0,0)[b]{\scriptsize $i$}}
\put(225,35){\makebox(0,0)[b]{\scriptsize $n\mn 2$}}

\put(150,29){\line(0,-1){18}} \put(180,29){\line(0,-1){18}}
\put(209.3,10.7){\line(-3,4){14}}
\put(239.3,10.7){\line(-3,4){14}}

\put(165,20){\makebox(0,0){\ldots}}
\put(220,20){\makebox(0,0){\ldots}}

\end{picture}
\end{center}
However, ${\cal H}\cirk{\cal J}$ is not faithful. For example,
$d_0^1\!:1\str 0$ and $d_1^1\!:1\str 0$ are both mapped to the
empty partial function from 1 to 0. So, $\cal H$ cannot be used
for the constructions like, for example, the functor nerve is,
where $d_0^1$ and $d_1^1$ should be mapped to the source and the
target function, respectively. The image of $\Delta_+^{op}$ in
$\Delta_{par}$ under ${\cal H}\cirk{\cal J}$ is a category
opposite to the category obtained as the image of $\Delta$ in the
category $\Gamma$ under a functor defined in \cite{S74}
(Section~1, first paragraph after Definition 1.2).

Let now $\langle C,\mu',\eta'\rangle$ be a monoid in a category
$\cal K$ whose monoidal structure is given by finite products. (We
abbreviate the product $C\times C$ by $C^2$, etc.) To avoid
permanent decoration with associativity, and right and left
identity isomorphisms of the monoidal structure of $\cal K$, we
will always consider this monoidal structure to be strict, which
is supported by the strictification given by \cite{ML71}, XI.3,
Theorem~1. Let $I$ be a terminal object of $\cal K$, which is the
strict monoidal unit. Then for $\varepsilon'$ being the unique
arrow from $C$ to $I$, we have that $\langle
C,\mu',\eta',\varepsilon'\rangle$ is a counital monoid. By
Proposition \thesection.2, we have a strict monoidal functor
${\cal F}\!:\Delta_{par}\str{\cal K}$ such that ${\cal F}(1)=C$,
etc.

If we denote by $\overline{W}C$ the composition ${\cal
F}\cirk{\cal H}\cirk{\cal J}\!:\Delta_+^{op}\str{\cal K}$, then we
have:
\begin{tabbing}
\hspace{1.5em}\=$\overline{W}C(d_0^n)$ \=$=pr_2\!:C\times
C^{n-1}\str C^{n-1}$,\hspace{1em}\=$\overline{W}C(d_n^n)$
\=$=pr_1\!:C^{n-1}\times C\str C^{n-1}$,
\\[2ex]
and for $1\leq i\leq n\mn 1$ and $0\leq j\leq n\mn 1$,
\\[1ex]
\>$\overline{W}C(d_i^n)$\>$=\mj^{i-1}\times\mu\times
\mj^{n-i-1}$,\>$\overline{W}C(s_j^n)$\>$=
\mj^j\times\eta\times\mj^{n-j-1}$.
\end{tabbing}
Hence, $\overline{W}C$ is the reduced bar construction of
\cite{T79}. Note that one needs just a part of the cartesian
structure of $\cal K$ that provides ``counits'', for such a
reduced bar construction. By relying on Proposition \thesection.2,
the construction may work in any strict monoidal category $\cal K$
equipped with a counital monoid $\langle
C,\mu',\eta',\varepsilon'\rangle$. Such a more general
construction of the composition ${\cal F}\cirk{\cal H}\cirk{\cal
J}\!:\Delta_+^{op}\str{\cal K}$ corresponds to ${\cal B}(I,C,I)$,
a special case of the two-sided bar construction in which we
regard $\cal K$ as a 2-category with one 0-cell, and $\langle
I,\varepsilon'\rangle$ as a right and left module over $C$.

However, for our purposes it is sufficient to consider the case
when $\cal K$ is the category \emph{Cat} (regarded again as strict
monoidal) and when $C$ is a strict monoidal category $\cal C$,
hence a monoid in \emph{Cat}. For example, if we take the arrow
$f\!:4\str 3$ of $\Delta_+^{op}$ graphically presented as

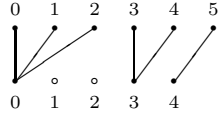
\begin{center}
\begin{picture}(85,40)

\put(0,10){\circle*{2}} \put(15,10){\circle{2}}
\put(30,10){\circle{2}} \put(45,10){\circle*{2}}
\put(60,10){\circle*{2}} \put(0,30){\circle*{2}}
\put(15,30){\circle*{2}} \put(30,30){\circle*{2}}
\put(45,30){\circle*{2}} \put(60,30){\circle*{2}}
\put(75,30){\circle*{2}}

\put(0,0){\makebox(0,0)[b]{\scriptsize $0$}}
\put(15,0){\makebox(0,0)[b]{\scriptsize $1$}}
\put(30,0){\makebox(0,0)[b]{\scriptsize $2$}}
\put(45,0){\makebox(0,0)[b]{\scriptsize $3$}}
\put(60,0){\makebox(0,0)[b]{\scriptsize $4$}}

\put(0,35){\makebox(0,0)[b]{\scriptsize $0$}}
\put(15,35){\makebox(0,0)[b]{\scriptsize $1$}}
\put(30,35){\makebox(0,0)[b]{\scriptsize $2$}}
\put(45,35){\makebox(0,0)[b]{\scriptsize $3$}}
\put(60,35){\makebox(0,0)[b]{\scriptsize $4$}}
\put(75,35){\makebox(0,0)[b]{\scriptsize $5$}}

\put(0,29){\line(0,-1){18}} \put(45,29){\line(0,-1){18}}

\put(0.7,10.7){\line(3,4){14}} \put(0.7,10.7){\line(3,2){29}}
\put(45.7,10.7){\line(3,4){14}} \put(60.7,10.7){\line(3,4){14}}

\end{picture}
\end{center}
then for $\overline{W}{\cal
C}\!:\Delta_+^{op}\str\mbox{\emph{Cat}}$, we have that
$\overline{W}{\cal C}(f)$, denoted by $f^\ast\!:{\cal
C}^4\str{\cal C}^3$, is a functor such that
\[
f^\ast(A,B,C,D)=(I,I,C\otimes D),
\]
where $I$ is the unit and $\otimes$ is the tensor of the strict
monoidal category $\cal C$.

\section{Iterated reduced bar construction}

Let now $\cal C$ be an $SMI$ category which is strict monoidal
with respect to both $\vee$, $\bot$ and $\kon$, $\top$. For $n$,
$m$ such that $n\pl m\geq 1$ we define, along the lines of
\cite{BFSV}, a lax functor (cf.\ \cite{S72}), which is an ordinary
functor when $n\pl m=1$,

\begin{center}
\begin{picture}(160,25)

\put(0,15){\makebox(0,0){$\overline{W}{\cal C}_{n,m}\!:$}}
\put(70,15){\makebox(0,0){$\Delta_+^{op}\times\ldots\times\Delta_+^{op}$}}
\put(70,9){\makebox(0,0)[b]{$\underbrace{\hspace{70pt}}$}}
\put(71,0){\makebox(0,0)[t]{\scriptsize $n\pl m$}}
\put(135,15){\makebox(0,0){$\str$}}
\put(160,15){\makebox(0,0){\emph{Cat}.}}

\end{picture}
\end{center}

It is defined on objects as $\overline{W}{\cal
C}_{n,m}(k_1,\ldots,k_{n+m})=_{df}{\cal C}^{k_1\cdot\ldots\cdot
k_{n+m}}$, and for the arrows we have the following. First, for
$1\leq i\leq n\pl m$ and $f_i\!:k_i\str l_i$ an arrow of
$\Delta_+^{op}$ let
\[
\overline{W}{\cal
C}_{n,m}(\mj_{k_1},\ldots,\mj_{k_{i-1}},f_i,\mj_{k_{i+1}},\ldots,\mj_{k_{n+m}})=_{df}
[\overline{W}{\cal D}(f_i)]^{k_1\cdot\ldots\cdot k_{i-1}},
\]
where $\cal D$ is the category ${\cal C}^{k_{i+1}\cdot\ldots\cdot
k_{n+m}}$ whose monoidal structure is defined componentwise in
terms of $\vee,\bot$ when $i\leq n$ and in terms of $\kon,\top$
when $n<i\leq n\pl m$. With this in mind, we define
$\overline{W}{\cal C}_{n,m}(f_1,\ldots,f_{n+m})$ to be the
following composition:
\[
\overline{W}{\cal C}_{n,m}(\mj_{l_1},\ldots,f_{n+m})
\cirk\ldots\cirk\overline{W}{\cal
C}_{n,m}(f_1,\mj_{k_2},\ldots,\mj_{k_{n+m}}).
\]

We call this construction of $\overline{W}{\cal C}_{n,m}$, the
$(n,m)$-\emph{reduced bar construction based on} $\cal C$. (When
$n\pl m=0$ we may define it to be the functor mapping the object
and the arrow of the trivial category $(\Delta_+^{op})^0$ to $\cal
C$ and to the identity functor on $\cal C$, respectively.)

\vspace{2ex} \noindent{\sc Example.} Let $n=2$, $m=1$, and let
$f=(f_1,f_2,f_3)\!:(1,2,2)\str(2,1,2)$ be the arrow of
$\Delta_+^{op}\times\Delta_+^{op}\times\Delta_+^{op}$ graphically
presented as

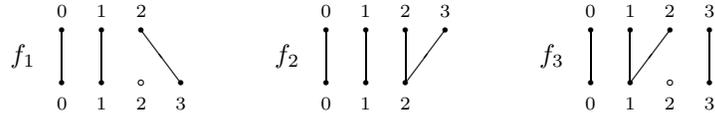
\begin{center}
\begin{picture}(245,40)

\put(0,10){\circle*{2}} \put(15,10){\circle*{2}}
\put(30,10){\circle{2}} \put(45,10){\circle*{2}}
\put(0,30){\circle*{2}} \put(15,30){\circle*{2}}
\put(30,30){\circle*{2}}

\put(0,0){\makebox(0,0)[b]{\scriptsize $0$}}
\put(15,0){\makebox(0,0)[b]{\scriptsize $1$}}
\put(30,0){\makebox(0,0)[b]{\scriptsize $2$}}
\put(45,0){\makebox(0,0)[b]{\scriptsize $3$}}

\put(0,35){\makebox(0,0)[b]{\scriptsize $0$}}
\put(15,35){\makebox(0,0)[b]{\scriptsize $1$}}
\put(30,35){\makebox(0,0)[b]{\scriptsize $2$}}

\put(-10,20){\makebox(0,0)[r]{$f_1$}}

\put(0,29){\line(0,-1){18}} \put(15,29){\line(0,-1){18}}

\put(44.3,10.7){\line(-3,4){14}}

\put(100,10){\circle*{2}} \put(115,10){\circle*{2}}
\put(130,10){\circle*{2}} \put(145,30){\circle*{2}}
\put(100,30){\circle*{2}} \put(115,30){\circle*{2}}
\put(130,30){\circle*{2}}

\put(100,0){\makebox(0,0)[b]{\scriptsize $0$}}
\put(115,0){\makebox(0,0)[b]{\scriptsize $1$}}
\put(130,0){\makebox(0,0)[b]{\scriptsize $2$}}
\put(145,35){\makebox(0,0)[b]{\scriptsize $3$}}

\put(100,35){\makebox(0,0)[b]{\scriptsize $0$}}
\put(115,35){\makebox(0,0)[b]{\scriptsize $1$}}
\put(130,35){\makebox(0,0)[b]{\scriptsize $2$}}

\put(90,20){\makebox(0,0)[r]{$f_2$}}

\put(100,29){\line(0,-1){18}} \put(115,29){\line(0,-1){18}}
\put(130,29){\line(0,-1){18}}

\put(130.7,10.7){\line(3,4){14}}

\put(200,10){\circle*{2}} \put(215,10){\circle*{2}}
\put(230,10){\circle{2}} \put(245,10){\circle*{2}}
\put(245,30){\circle*{2}} \put(200,30){\circle*{2}}
\put(215,30){\circle*{2}} \put(230,30){\circle*{2}}

\put(200,0){\makebox(0,0)[b]{\scriptsize $0$}}
\put(215,0){\makebox(0,0)[b]{\scriptsize $1$}}
\put(230,0){\makebox(0,0)[b]{\scriptsize $2$}}
\put(245,0){\makebox(0,0)[b]{\scriptsize $3$}}

\put(200,35){\makebox(0,0)[b]{\scriptsize $0$}}
\put(215,35){\makebox(0,0)[b]{\scriptsize $1$}}
\put(230,35){\makebox(0,0)[b]{\scriptsize $2$}}
\put(245,35){\makebox(0,0)[b]{\scriptsize $3$}}

\put(190,20){\makebox(0,0)[r]{$f_3$}}

\put(200,29){\line(0,-1){18}} \put(245,29){\line(0,-1){18}}
\put(215,29){\line(0,-1){18}}

\put(215.7,10.7){\line(3,4){14}}

\end{picture}
\end{center}
and $g=(g_1,g_2,g_3)\!:(2,1,2)\str(2,2,2)$ be the arrow of
$\Delta_+^{op}\times\Delta_+^{op}\times\Delta_+^{op}$ graphically
presented as

\begin{center}
\begin{picture}(245,40)

\put(0,10){\circle*{2}} \put(15,10){\circle*{2}}
\put(30,10){\circle{2}} \put(45,10){\circle*{2}}
\put(45,30){\circle*{2}} \put(0,30){\circle*{2}}
\put(15,30){\circle*{2}} \put(30,30){\circle*{2}}

\put(0,0){\makebox(0,0)[b]{\scriptsize $0$}}
\put(15,0){\makebox(0,0)[b]{\scriptsize $1$}}
\put(30,0){\makebox(0,0)[b]{\scriptsize $2$}}
\put(45,0){\makebox(0,0)[b]{\scriptsize $3$}}

\put(0,35){\makebox(0,0)[b]{\scriptsize $0$}}
\put(15,35){\makebox(0,0)[b]{\scriptsize $1$}}
\put(30,35){\makebox(0,0)[b]{\scriptsize $2$}}
\put(45,35){\makebox(0,0)[b]{\scriptsize $3$}}

\put(-10,20){\makebox(0,0)[r]{$g_1$}}

\put(0,29){\line(0,-1){18}} \put(45,29){\line(0,-1){18}}
\put(15,29){\line(0,-1){18}}

\put(15.7,10.7){\line(3,4){14}}

\put(100,10){\circle*{2}} \put(115,10){\circle*{2}}
\put(130,10){\circle{2}} \put(145,10){\circle*{2}}
\put(100,30){\circle*{2}} \put(115,30){\circle*{2}}
\put(130,30){\circle*{2}}

\put(100,0){\makebox(0,0)[b]{\scriptsize $0$}}
\put(115,0){\makebox(0,0)[b]{\scriptsize $1$}}
\put(130,0){\makebox(0,0)[b]{\scriptsize $2$}}
\put(145,0){\makebox(0,0)[b]{\scriptsize $3$}}

\put(100,35){\makebox(0,0)[b]{\scriptsize $0$}}
\put(115,35){\makebox(0,0)[b]{\scriptsize $1$}}
\put(130,35){\makebox(0,0)[b]{\scriptsize $2$}}

\put(90,20){\makebox(0,0)[r]{$g_2$}}

\put(100,29){\line(0,-1){18}} \put(115,29){\line(0,-1){18}}

\put(144.3,10.7){\line(-3,4){14}}

\put(200,10){\circle*{2}} \put(215,10){\circle*{2}}
\put(230,10){\circle*{2}} \put(245,10){\circle*{2}}
\put(245,30){\circle*{2}} \put(200,30){\circle*{2}}
\put(215,30){\circle*{2}} \put(230,30){\circle*{2}}

\put(200,0){\makebox(0,0)[b]{\scriptsize $0$}}
\put(215,0){\makebox(0,0)[b]{\scriptsize $1$}}
\put(230,0){\makebox(0,0)[b]{\scriptsize $2$}}
\put(245,0){\makebox(0,0)[b]{\scriptsize $3$}}

\put(200,35){\makebox(0,0)[b]{\scriptsize $0$}}
\put(215,35){\makebox(0,0)[b]{\scriptsize $1$}}
\put(230,35){\makebox(0,0)[b]{\scriptsize $2$}}
\put(245,35){\makebox(0,0)[b]{\scriptsize $3$}}

\put(190,20){\makebox(0,0)[r]{$g_3$}}

\put(200,29){\line(0,-1){18}} \put(230,29){\line(0,-1){18}}
\put(215,29){\line(0,-1){18}} \put(245,29){\line(0,-1){18}}

\end{picture}
\end{center}

With abbreviation $h^\ast$ for $\overline{W}{\cal C}_{2,1}(h)$, we
have that $f^\ast$ is the functor from ${\cal C}^{1\cdot 2\cdot
2}$ to ${\cal C}^{2\cdot 1\cdot 2}$ defined by
\[
(A,B,C,D)\mapsto(A,B,C,D,\bot,\bot,\bot,\bot)\mapsto(A,B,\bot,\bot)\mapsto(A\kon
B,\top,\bot\kon\bot,\top),
\]
and $g^\ast$ is the functor from ${\cal C}^{2\cdot 1\cdot 2}$ to
${\cal C}^{2\cdot 2\cdot 2}$ defined by
\[
(A,B,C,D)\mapsto(A\vee C,B\vee D,\bot,\bot)\mapsto(A\vee C,B\vee
D,\bot,\bot,\bot,\bot,\bot,\bot).
\]

That $\overline{W}{\cal C}_{2,1}$ is not a functor could be seen
from the fact that $g^\ast\cirk f^\ast$, which is defined by
\[
(A,B,C,D)\mapsto((A\kon
B)\vee(\bot\kon\bot),\top\vee\top,\bot,\bot,\bot,\bot,\bot,\bot),
\]
is different from $(g\cirk f)^\ast$. Here $g\cirk f=(g_1\cirk
f_1,g_2\cirk f_2,g_3\cirk f_3)\!:(1,2,2)\str(2,2,2)$ is
graphically presented as

\begin{center}
\begin{picture}(245,40)

\put(0,10){\circle*{2}} \put(15,10){\circle*{2}}
\put(30,10){\circle{2}} \put(45,10){\circle*{2}}
\put(0,30){\circle*{2}} \put(15,30){\circle*{2}}
\put(30,30){\circle*{2}}

\put(0,0){\makebox(0,0)[b]{\scriptsize $0$}}
\put(15,0){\makebox(0,0)[b]{\scriptsize $1$}}
\put(30,0){\makebox(0,0)[b]{\scriptsize $2$}}
\put(45,0){\makebox(0,0)[b]{\scriptsize $3$}}

\put(0,35){\makebox(0,0)[b]{\scriptsize $0$}}
\put(15,35){\makebox(0,0)[b]{\scriptsize $1$}}
\put(30,35){\makebox(0,0)[b]{\scriptsize $2$}}

\put(-7,20){\makebox(0,0)[r]{$g_1\cirk f_1$}}

\put(0,29){\line(0,-1){18}} \put(15,29){\line(0,-1){18}}

\put(44.3,10.7){\line(-3,4){14}} \put(100,10){\circle*{2}}
\put(115,10){\circle*{2}} \put(130,10){\circle{2}}
\put(145,10){\circle*{2}} \put(100,30){\circle*{2}}
\put(115,30){\circle*{2}} \put(130,30){\circle*{2}}
\put(145,30){\circle*{2}}

\put(100,0){\makebox(0,0)[b]{\scriptsize $0$}}
\put(115,0){\makebox(0,0)[b]{\scriptsize $1$}}
\put(130,0){\makebox(0,0)[b]{\scriptsize $2$}}
\put(145,0){\makebox(0,0)[b]{\scriptsize $3$}}

\put(100,35){\makebox(0,0)[b]{\scriptsize $0$}}
\put(115,35){\makebox(0,0)[b]{\scriptsize $1$}}
\put(130,35){\makebox(0,0)[b]{\scriptsize $2$}}
\put(145,35){\makebox(0,0)[b]{\scriptsize $3$}}

\put(93,20){\makebox(0,0)[r]{$g_2\cirk f_2$}}

\put(100,29){\line(0,-1){18}} \put(115,29){\line(0,-1){18}}
\put(145,29){\line(0,-1){18}}

\put(144.3,10.7){\line(-3,4){14}}

\put(200,10){\circle*{2}} \put(215,10){\circle*{2}}
\put(230,10){\circle{2}} \put(245,10){\circle*{2}}
\put(245,30){\circle*{2}} \put(200,30){\circle*{2}}
\put(215,30){\circle*{2}} \put(230,30){\circle*{2}}

\put(200,0){\makebox(0,0)[b]{\scriptsize $0$}}
\put(215,0){\makebox(0,0)[b]{\scriptsize $1$}}
\put(230,0){\makebox(0,0)[b]{\scriptsize $2$}}
\put(245,0){\makebox(0,0)[b]{\scriptsize $3$}}

\put(200,35){\makebox(0,0)[b]{\scriptsize $0$}}
\put(215,35){\makebox(0,0)[b]{\scriptsize $1$}}
\put(230,35){\makebox(0,0)[b]{\scriptsize $2$}}
\put(245,35){\makebox(0,0)[b]{\scriptsize $3$}}

\put(193,20){\makebox(0,0)[r]{$g_3\cirk f_3$}}

\put(200,29){\line(0,-1){18}} \put(245,29){\line(0,-1){18}}
\put(215,29){\line(0,-1){18}}

\put(215.7,10.7){\line(3,4){14}}

\end{picture}
\end{center}
and $(g\cirk f)^\ast\!:{\cal C}^{1\cdot 2\cdot 2}\str{\cal
C}^{2\cdot 2\cdot 2}$ is defined by
\begin{tabbing}
\hspace{2em}$(A,B,C,D)\mapsto(A,B,C,D,\bot,\bot,\bot,\bot)\mapsto
(A,B,\bot,\bot,\bot,\bot,\bot,\bot)\mapsto$
\\[2ex]
\`$\mapsto(A\kon B,\top,\bot\kon\bot,
\top,\bot\kon\bot,\top,\bot\kon\bot,\top)$.
\end{tabbing}

However, we have a natural transformation from $g^\ast\cirk
f^\ast$ to $(g\cirk f)^\ast$ whose components are
$\iota_{A,B,\bot,\bot}$, $\tau^\str$, $\beta^\rts$ (three times)
and $\kappa$ (three times). This natural transformation acts as
$\omega_{g,f}$ from the definition of lax functor and since
$\overline{W}{\cal C}_{n,m}$ preserves identity arrows, there is
no need for the natural transformation $\omega_A$.

To show that $\overline{W}{\cal C}_{n,m}$, for $n\pl m\geq 2$, is
indeed a lax functor, we have to find for every composable pair of
arrows $f$, $g$ of $\Delta_+^{op}\times\ldots\times\Delta_+^{op}$,
a natural transformation $\omega_{g,f}\!:g^\ast\cirk
f^\ast\strt(g\cirk f)^\ast$, such that the following diagram
commutes

\begin{center}
\begin{picture}(100,80)

\put(-50,40){\makebox(0,0)[r]{(\emph{lax})}}
\put(50,70){\makebox(0,0){$h^\ast\cirk g^\ast\cirk f^\ast$}}
\put(0,40){\makebox(0,0){$(h\cirk g)^\ast\cirk f^\ast$}}
\put(100,40){\makebox(0,0){$h^\ast\cirk (g\cirk f)^\ast$}}
\put(50,10){\makebox(0,0){$(h\cirk g\cirk f)^\ast$}}

\put(30,60){\vector(-2,-1){20}} \put(70,60){\vector(2,-1){20}}
\put(10,30){\vector(2,-1){20}} \put(90,30){\vector(-2,-1){20}}

\put(20,60){\makebox(0,0)[r]{\scriptsize $\omega_{h,g}(f^\ast)$}}
\put(80,60){\makebox(0,0)[l]{\scriptsize $h^\ast(\omega_{g,f})$}}

\put(20,20){\makebox(0,0)[r]{\scriptsize $\omega_{h\cirk g,f}$}}
\put(80,20){\makebox(0,0)[l]{\scriptsize $\omega_{h,g\cirk f}$}}

\end{picture}
\end{center}
For this, we rely on the category \cmst, which is strict monoidal
$SMI$ category freely generated by the same infinite set $\cal P$
of generators we have used for the category \cm\ in Section~3. The
category \cmst\ is obtained from our category \cm\ by factoring
its objects through the smallest equivalence relation $\equiv$
satisfying
\begin{tabbing}
\hspace{2em}$A\vee(B\vee C)\equiv(A\vee B)\vee
C$,\hspace{5em}$A\kon(B\kon C)\equiv(A\kon B)\kon C$,
\\[1ex]
\centerline{$A\equiv A\kon\top\equiv\top\kon A\equiv
A\vee\bot\equiv\bot\vee A$,}
\end{tabbing}
which is congruent with respect to $\vee$ and $\kon$, and by
further factoring its arrow terms according to the new equations
\begin{tabbing}
\hspace{2em}$\check{\alpha}^\str_{A,B,C}=\check{\alpha}^\rts_{A,B,C}=\mj_{A\vee
B\vee
C}$,\hspace{5em}$\hat{\alpha}^\str_{A,B,C}=\hat{\alpha}^\rts_{A,B,C}=\mj_{A\kon
B\kon C}$,
\\[2ex]
\centerline{$\check{\rho}^\str_A=\check{\rho}^\rts_A=
\check{\lambda}^\str_A=\check{\lambda}^\rts_A=
\hat{\rho}^\str_A=\hat{\rho}^\rts_A=
\hat{\lambda}^\str_A=\hat{\lambda}^\rts_A=\mj_A$.}
\end{tabbing}
(Hence, while writing down the objects of \cmst, we may omit
parentheses tied to $\vee$, and the constant $\bot$ in the
immediate scope of another $\vee$ and the same for $\kon$ and
$\top$.) An object of \cmst\ is \emph{pure} and \emph{diversified}
when it, as an equivalence class, consists of formulae that are
pure and diversified. As a direct consequence of Theorem 3.1 we
have

\prop{Corollary of Theorem 3.1}{If $A$ and $B$ are either pure and
diversified or no letter occurs in them, then there is at most one
arrow $f\!:A\str B$ in \cmst.}

The following lemma reduces our problem to the category \cmst. The
hypotheses (1) in the formulation serves for the existence of an
arrow $\omega_{f,g}(\vec{A})$, whereas (2) together with the
coherence theorem are used to guarantee both its uniqueness and
the commutativity of the diagram (\emph{lax}) above.

\prop{Lemma \thesection.1}{If for $\overline{W}\cm^{st}_{n,m}$ the following holds:\\[1ex] \mbox{\rm (1)}\quad for every pair of arrows
$f\!:(k_1,\ldots,k_{n+m})\str(l_1,\ldots,l_{n+m})$ and \linebreak
$g\!:(l_1,\ldots,l_{n+m})\str (j_1,\ldots,j_{n+m})$ of
$\Delta_+^{op}\times\ldots\times\Delta_+^{op}$, and every
$k_1\cdot\ldots\cdot k_{n+m}$-tuple of different letters
$\vec{p}=p_{11\ldots1},p_{11\ldots2},\ldots,p_{k_1k_2\ldots
k_{n+m}}$, there is an arrow $\omega_{f,g}(\vec{p})\!:g^\ast\cirk
f^\ast(\vec{p})\str(g\cirk f)^\ast(\vec{p})$ of
$(\cm^{st})^{j_1\cdot\ldots\cdot j_{n+m}}$,
and \\[1ex] \mbox{\rm (2)}\quad for every sequence of composable arrows
$f_1\ldots f_u$ of $\Delta_+^{op}\times\ldots\times\Delta_+^{op}$,
each coordinate of $f_u^\ast\cirk\ldots\cirk f_1^\ast(\vec{p})$,
is
either pure and diversified or no letter occurs in it,\\[1ex] then
for every strict monoidal $SMI$ category $\cal C$, we have that
$\overline{W}{\cal C}_{n,m}$ is a lax functor. }

\dkz Using the freedom of \cmst\ and (1) we define for every
$k_1\cdot\ldots\cdot k_{n+m}$-tuple
$\vec{A}=(A_{11\ldots1},A_{11\ldots2},\ldots,A_{k_1\ldots
k_{n+m}})$ of objects of $\cal C$ the arrow
$\omega_{f,g}(\vec{A})$ as the image of $\omega_{f,g}(\vec{p})$
under the functor that extends the function mapping a generator
$p_{i_1\ldots i_{n+m}}$ of \cmst\ to the object $A_{i_1\ldots
i_{n+m}}$ of $\cal C$. From (2) ($u=1$ and $u=3$ are the only
interesting cases), appealing again to the freedom of \cmst, and
to Corollary of Theorem 3.1, we have that the diagram (\emph{lax})
commutes. \qed

\vspace{2ex}

To prove that (1) holds, we reason as in \cite{BFSV}. Since for
every $i\in\{1,\ldots,n\pl m\}$ we have that
\[
(\mj_{k_1},\ldots,\mj_{k_{i-1}},f_i,\mj_{k_{i+1}},\ldots,\mj_{k_{n+m}})^\ast
\]
is a functor, it is sufficient to show that for every $1\leq
i<j\leq n+m$, and $f_i\!:k_i\str l_i$ and  $f_j\!:k_j\str l_j$
arrows of $\Delta_+^{op}$ there is an arrow of
$(\cm^{st})^{k_1\cdot\ldots\cdot l_i\cdot\ldots\cdot
l_j\cdot\ldots\cdot k_{n+m}}$ whose source is
\[
(\mj_{k_1},\ldots,f_i,\ldots,\mj_{l_j},\ldots,\mj_{k_{n+m}})^\ast\cirk
(\mj_{k_1},\ldots,\mj_{k_i},\ldots,f_j,\ldots,\mj_{k_{n+m}})^\ast
(\vec{p}),
\]
and whose target is
\[
(\mj_{k_1},\ldots,\mj_{l_i},\ldots,f_j,\ldots,\mj_{k_{n+m}})^\ast\cirk
(\mj_{k_1},\ldots,f_i,\ldots,\mj_{k_j},\ldots,\mj_{k_{n+m}})^\ast
(\vec{p}).
\]
Since it is sufficient to find each coordinate of this arrow, we
may assume that all numbers except $k_i$ and $k_j$ are 1, and we
write $k_i\cdot k_j$ tuple $\vec{p}$ as
$p_{11},p_{12},\ldots,p_{k_ik_j}$.

Let $f_i\!:k_i\str 1$ and $f_j\!:k_j\str 1$ be the arrows of
$\Delta_+^{op}$ graphically presented as

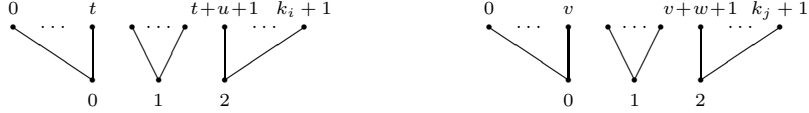
\begin{center}
\begin{picture}(290,40)

\put(30,10){\circle*{2}} \put(55,10){\circle*{2}}
\put(80,10){\circle*{2}} \put(0,30){\circle*{2}}
\put(30,30){\circle*{2}} \put(45,30){\circle*{2}}
\put(65,30){\circle*{2}} \put(80,30){\circle*{2}}
\put(110,30){\circle*{2}}

\put(30,0){\makebox(0,0)[b]{\scriptsize $0$}}
\put(55,0){\makebox(0,0)[b]{\scriptsize $1$}}
\put(80,0){\makebox(0,0)[b]{\scriptsize $2$}}

\put(0,35){\makebox(0,0)[b]{\scriptsize $0$}}
\put(30,35){\makebox(0,0)[b]{\scriptsize $t$}}
\put(80,35){\makebox(0,0)[b]{\scriptsize $t\pl u\pl 1$}}
\put(110,35){\makebox(0,0)[b]{\scriptsize $k_i+1$}}

\put(15,30){\makebox(0,0){\scriptsize $\ldots$}}
\put(55,30){\makebox(0,0){\scriptsize $\ldots$}}
\put(95,30){\makebox(0,0){\scriptsize $\ldots$}}

\put(30,29){\line(0,-1){18}} \put(80,29){\line(0,-1){18}}
\put(0,30){\line(3,-2){30}} \put(110,30){\line(-3,-2){30}}
\put(45,30){\line(1,-2){10}} \put(65,30){\line(-1,-2){10}}

\put(210,10){\circle*{2}} \put(235,10){\circle*{2}}
\put(260,10){\circle*{2}} \put(180,30){\circle*{2}}
\put(210,30){\circle*{2}} \put(225,30){\circle*{2}}
\put(245,30){\circle*{2}} \put(260,30){\circle*{2}}
\put(290,30){\circle*{2}}

\put(210,0){\makebox(0,0)[b]{\scriptsize $0$}}
\put(235,0){\makebox(0,0)[b]{\scriptsize $1$}}
\put(260,0){\makebox(0,0)[b]{\scriptsize $2$}}

\put(180,35){\makebox(0,0)[b]{\scriptsize $0$}}
\put(210,35){\makebox(0,0)[b]{\scriptsize $v$}}
\put(260,35){\makebox(0,0)[b]{\scriptsize $v\pl w\pl 1$}}
\put(290,35){\makebox(0,0)[b]{\scriptsize $k_j+1$}}

\put(195,30){\makebox(0,0){\scriptsize $\ldots$}}
\put(235,30){\makebox(0,0){\scriptsize $\ldots$}}
\put(275,30){\makebox(0,0){\scriptsize $\ldots$}}

\put(210,29){\line(0,-1){18}} \put(260,29){\line(0,-1){18}}
\put(180,30){\line(3,-2){30}} \put(290,30){\line(-3,-2){30}}
\put(225,30){\line(1,-2){10}} \put(245,30){\line(-1,-2){10}}

\end{picture}
\end{center}
So, in the case when $i<j\leq n$ we need an arrow
\[
\bigvee_{x=v+1}^{v+w} \bigvee_{y=t+1}^{t+u}p_{xy}\str
\bigvee_{y=t+1}^{t+u}\bigvee_{x=v+1}^{v+w}p_{xy},
\]
which is $\mj_\bot$ when either $u$ or $w$ is 0, or it is built
out of $\check{\sigma}$, otherwise. In the case when $i\leq n<j$
we need an arrow
\[
\bigvee_{x=v+1}^{v+w}\bigwedge_{y=t+1}^{t+u}p_{xy}\str
\bigwedge_{y=t+1}^{t+u}\bigvee_{x=v+1}^{v+w}p_{xy},
\]
which is built out of $\tau^\str$, $\beta^\rts$ and $\kappa$ when
$u$ or $v$ is 0, or it is built out of $\iota$, otherwise. In the
case when $n<i<j$ we proceed as in the first case relying on
$\mj_\top$ and $\hat{\sigma}$. So, (1) is proved.

To prove that (2) holds, note first that the equivalence relation
used to factor the objects of \cm\ in order to obtain the objects
of \cmst\ is congruent with respect to the function $\nu$ defined
in Section~3. So, $\nu$ may be considered as a function on the
objects of \cmst. We say that an object
$\vec{A}=(A_{11\ldots1},A_{11\ldots2},\ldots,A_{k_1k_2\ldots
k_{n+m}})$ of $(\cm^{st})^{k_1\cdot\ldots\cdot k_{n+m}}$ is
$(n,m)$-\emph{coherent} when the following holds for $1\leq
i_l,j_l\leq k_l$:

\vspace{1ex}

$(\ast)$ Every $A_{i_1i_2\ldots i_{n+m}}$ is either pure and
diversified or no letter occurs in it, and $let(A_{i_1i_2\ldots
i_{n+m}})\cap let(A_{j_1j_2\ldots j_{n+m}})=\emptyset$, when
$i_1i_2\ldots i_{n+m}\neq j_1j_2\ldots j_{n+m}$;

\vspace{1ex}

$(\ast\ast)$ For every $m$-tuple $i_{n+1}\ldots i_{n+m}$, if for
some $n$-tuple $i_1,\ldots,i_n$ we have that $\nu(A_{i_1\ldots i_n
i_{n+1}\ldots i_{n+m}})$ is $\top$, then for every $n$-tuple
$j_1,\ldots,j_n$ we have that $\nu(A_{j_1\ldots j_n i_{n+1}\ldots
i_{n+m}})$ is $\top$ or $\bot$;

\vspace{1ex}

$(\ast\ast\ast)$ For every $n$-tuple $i_1,\ldots,i_n$, if for some
$m$-tuple $i_{n+1},\ldots,i_{n+m}$ we have that $\nu(A_{i_1\ldots
i_n i_{n+1}\ldots i_{n+m}})$ is $\bot$, then for every $m$-tuple
$j_{n+1},\ldots,j_{n+m}$ we have that $\nu(A_{i_1\ldots i_n
j_{n+1}\ldots j_{n+m}})$ is $\top$ or $\bot$.

\vspace{1ex}

The following lemma has (2) as an immediate corollary.

\prop{Lemma \thesection.2}{For every $1\leq i\leq n\pl m$ and
every arrow $f_i\!:k_i\str l_i$ of $\Delta_+^{op}$, if $\vec{A}$
is an $(n,m)$-coherent object of $(\cm^{st})^{k_1\cdot\ldots\cdot
k_{n+m}}$, then
\[ \vec{B}=(\mj_{k_1},\ldots,\mj_{k_{i-1}},f_i,\ldots,\mj_{k_{n+m}})^\ast
(\vec{A})\] is an $(n,m)$-coherent object of
$(\cm^{st})^{k_1\cdot\ldots\cdot l_i\cdot\ldots\cdot k_{n+m}}$.}

\dkz Since for every $i\in\{1,\ldots,n\pl m\}$ we have that
\[
(\mj_{k_1},\ldots,\mj_{k_{i-1}},f_i,\mj_{k_{i+1}},\ldots,\mj_{k_{n+m}})^\ast
\]
is a functor, it is sufficient to prove the lemma for $f_i$ being
$d^{k_i}_j\!:k_i\str k_i\mn1$ or $s^{k_i+1}_j\!:k_i\str k_i\pl1$.
One can use the following table to verify that $\vec{B}$ satisfies
$(\ast)$, $(\ast\ast)$ and $(\ast\ast\ast)$. In this table the
index $\alpha$ is an $n\pl m$ sequence of natural numbers,
$\alpha_i\in[1,l_i]$ is its $i$-th component, $e_i$ is the $n\pl
m$ sequence with 1 as the $i$-th component and 0 everywhere else,
and the addition-subtraction is componentwise.

\begin{center}
\begin{tabular}{c|cl}
$f_i$ & $B_\alpha$ &
\\[.5ex]
\hline\hline\\[-2.2ex] $d_0^{k_i}$ & $A_{\alpha+e_i}$ &
\\[.5ex]
\hline\\[-2.2ex] $d_{k_i}^{k_i}$ & $A_\alpha$ &
\\[.5ex]
\hline\\[-2.2ex] & $A_\alpha$ & $\alpha_i<j$
\\
$d_j^{k_i}$ & $A_\alpha\vee A_{\alpha+e_i}$ & $\alpha_i=j$ \&
$1\leq i\leq n$\quad $(\ast\ast)$
\\ $0<j<k_i$ & $A_\alpha\kon A_{\alpha+e_i}$ & $\alpha_i=j$ \& $n<
i\leq n\pl m$\quad $(\ast\ast\ast)$
\\ & $A_{\alpha+e_i}$ & $\alpha_i>j$
\\[.5ex]
\hline\\[-2.2ex] & $A_\alpha$ & $\alpha_i<j\pl 1$
\\
$s_j^{k_i+1}$ & $\bot$ & $\alpha_i=j\pl 1$ \& $1\leq i\leq n$
\\ $0\leq j\leq k_i$ & $\top$ & $\alpha_i=j\pl 1$ \& $n<i\leq n\pl m$
\\ & $A_{\alpha-e_i}$ & $\alpha_i>j\pl 1$

\end{tabular}
\end{center}

In the case marked by $(\ast\ast)$ we use essentially the property
$(\ast\ast)$ of $\vec{A}$ to establish that $(\ast)$ holds for
$\vec{B}$, and analogously for $(\ast\ast\ast)$. This is the
reason why the properties $(\ast\ast)$ and $(\ast\ast\ast)$ occur
in the definition of an $(n,m)$-coherent
object.\mbox{\hspace{1em}} \qed

Now (2) follows immediately since every $k_1\cdot\ldots\cdot
k_{n+m}$-tuple $\vec{p}$ of different letters is obviously
$(n,m)$-coherent and one has just to iterate Lemma \thesection.2
through the definition of $f_u^\ast\cirk\ldots\cirk
f_1^\ast(\vec{p})$, and eventually, to use the property $(\ast)$
of the obtained object. Hence, we conclude from Lemma
\thesection.1 that every $(n,m)$-reduced bar construction based on
a strict monoidal $SMI$ category $\cal C$ produces a lax functor.

\section{Two questions}

The following questions, to which we have no answer, may come to
mind to a careful reader of this paper:

1) Do we have an unrestricted coherence for $SMI$ categories, i.e.
whether all diagrams (with diversified objects in the nodes)
commute in \cm?

2) Since there is no need for $\tau^\rts$ and $\beta^\str$ in the
construction of $\omega_{f,g}$, is it possible to omit the
assumptions that $\tau^\str$ and $\beta^\rts$ are isomorphisms
from the definition of $SMI$ categories, without loss of coherence
necessary for $(n,m)$-reduced bar construction?

The first question is of lower interest, at least for the
$(n,m)$-reduced bar construction, since we managed to work without
unrestricted coherence for $SMI$ categories. Some serious doubts
about holding of such a coherence result may be found in
\cite{DP08c}, Section~7. However, the second question may be quite
interesting for the matters of $(n,m)$-reduced bar construction.
An affirmative answer says that this construction may be based on
an arbitrary object of $SyMon_{lax}^{2}(Cat)$ (cf.\ Section~1),
and hence on every category with finite coproducts and products
without restriction to those categories having initial object as
the product of initial objects, and terminal object as the
coproduct of terminal objects.

\end{document}